\documentclass[11pt]{article}
\usepackage{amsmath, amssymb}
\usepackage[french]{babel}
\title{Character expansion method
for the first order asymptotics of  a matrix integral 
}
\author{Alice Guionnet\thanks{Ecole Normale Sup\'erieure de Lyon,
Unit\'e de Math\'ematiques pures et appliqu\'ees,
UMR 5669,
46 All\'ee d'Italie, 
69364 Lyon Cedex 07, France. E-mail: aguionne@umpa.ens-lyon.fr.}, Myl\`ene Ma\"{\i}da\thanks{Ecole Normale Sup\'erieure de Lyon,
Unit\'e de Math\'ematiques pures et appliqu\'ees,
UMR 5669,
46 All\'ee d'Italie, 
69364 Lyon Cedex 07, France. E-mail: mmaida@umpa.ens-lyon.fr}}

\newtheorem{prop}{Proposition}[section]
\newtheorem{theo}[prop]{Theorem}
\newtheorem{defi}[prop]{Definition}
\newtheorem{lem}[prop]{Lemma}

\newtheorem{hyp}[prop]{Hypothesis}
\newtheorem{rmk}[prop]{Remark}
\newcommand{\RR}{\mathbb{R}}
\newcommand{\CC}{\mathbb{C}}
\newcommand{\NN}{\mathbb{N}}

\newcommand{\PR}{\mathcal{P}(\RR)}
\newcommand{\Ppos}{\mathcal{P}(\RR^{+})}

\newcommand{\MNC}{\mathcal{M}_N(\CC)}
\newcommand{\HNC}{\mathcal{H}_N(\CC)}
\newcommand{\UNC}{\mathcal{U}_N(\CC)}

\newcommand{\ppq}{\leqslant}
\newcommand{\pgq}{\geqslant}

\def\prf{{\bf Proof.}}

\def\tr{{\mbox{tr}}}

\def\nn{\noindent}

\def\eps{\epsilon}

\def\xx{\vrule height 0.7em depth 0.2em width 0.5 em}

\def\a{\alpha}
\def\b{\beta}
\def\d{\delta}
\def\e{\epsilon}

\def\l{\lambda}

\def\D{\Delta}

\def\Si{\Sigma}

\def\ra{\rightarrow}

\def\Aa{{\cal A}}

\def\Ca{{\cal C}}

\def\Ga{{\cal G}}
\def\Ha{{\cal H}}
\def\Ia{{\cal I}}

\def\Ka{{\cal K}}

\def\La{{\cal L}}
\def\Ma{{\cal M}}

\def\Za{{\cal Z}}

\def\mun{{\hat\mu^N}}
\def\R{\RR}

\def\lbc{\lbrace}

\def\nn{\noindent}
\def\part{\partial}

\def\ot{\otimes}

\def\ts{\times}
\def\Pa{{\mathcal{P}}}
\begin{document}

\renewcommand{\refname}{\Large{References}}

\maketitle

\centerline{\bf Abstract}
The estimation of various matrix integrals 
as the size of the matrices goes to infinity 
is motivated  by theoretical
physics, geometry and free probability questions.
On a rigorous ground, only integrals of one matrix or
of several matrices with simple quadratic 
interaction (called $AB$ interaction)
could be evaluated so far (see e.g. \cite{M-M}, 
\cite{Me2} or \cite{G}).
In this article, we follow
an idea widely developed in
the physics literature, which
is based on character 
expansion, to study more complex interaction.
In this context, we derive a large deviation principle 
for the empirical measure of Young tableaux. We then use it
to  study a matrix  model
 defined in the spirit
of the 'dually weighted graph model' introduced 
in \cite{K2}, but with a cutoff function 
such that the matrix integral 
and its character expansion converge. We prove 
that the free energy of this model
converges as the size of
the matrices go to infinity  and study 
the saddle points of the 
limit.

\medskip

\nn
{\it Keywords :} Large deviations, random matrices, non-commutative measure,
integration.

\medskip

\nn
{\it Mathematics Subject of Classification :}  60F10, 15A52, 46L50.

\nn

\section{Introduction}
The evaluation of matrix integrals
was first motivated by theoretical
physics and geometry since 
they can be  related, via Feynman diagrams expansion (see \cite{ZV}
for a nice introduction), to
the enumeration of maps. Thanks to this relation,
matrix integrals can also be used
to describe some models appearing in statistical mechanics,
such as the Ising model or
the q-Potts model, on random graphs (instead
of the usual two-dimensional lattice). Using similar 
ideas, string theory 
models can be described via 
matrix integrals around criticality (see the course \cite{Ey} 
for various applications to physics).
Another motivation is the study of
non-commutative entropies introduced by
D. Voiculescu \cite{Vo1} in the context 
of free probability. Let us roughly say 
that the understanding
of the asymptotic behavior 
of all  possible matrix integrals 
would be equivalent to
the understanding of the so-called 
microstates entropy.\\
 So, what is a matrix integral ?
If we let, for $n\in\NN$, $\CC\langle X_1,\cdots, X_n\rangle$
be the set of  polynomial functions of $n$ 
non-commutative variables and if we choose, for some $m,p\in\NN$,
$P\in \CC\langle X_1,\cdots, X_{n+p}\rangle^{\ot m}$
and $\phi:=(\phi_i)_{1\le i\le n+p}\in\Ca^o(\R)^{n+p}$,
then a matrix integral can be defined by

$$Z_N(P,\phi)=\int e^{N^2 (N^{-1}\tr)^{\ot m}(P(\phi_1(A_1),\cdots, \phi_{n+p}(A_{n+p}))}
dA_1\cdots d A_n.$$
where $dA$ denotes the Lebesgue measure
on the chosen state space of the matrices, included into $\Ma_N
(\CC)$, the space of square matrices of dimension $N$ with complex entries.
In the following, the matrices will
take their values in the set $\Ha_N(\CC)$ 
of  Hermitian matrices of dimension $N$.
The first order asymptotics of $Z_N(P,\phi)$
can easily be studied in the case where $n=1$ 
since then the joint law of
the eigenvalues of the matrix $A$ 
is known and described by the Coulomb gas law (see
\cite{BA-G} for instance). All
the correction terms have been recently studied rigorously
by N. Ercolani and K. McLaughlin 
in \cite{EML}. To this end, they
use Riemann-Hilbert techniques together 
with a good understanding of the asymptotic behaviour
of the spectral measure of the matrix
with law given by the corresponding Gibbs measure
$$\mu_N^{P, \phi}(dA_1\cdots d A_n)=Z_N(P,\phi)^{-1}
e^{N^2 (N^{-1}\tr)^{\ot m}(P(\phi_1(A_1),\cdots, \phi_{n+p}(A_{n+p}))}
dA_1\cdots d A_n
.$$
There are much less complete results
in the case where $n\ge 2$. On a rigorous
ground,
let us however mention the work
of M. Mehta and al. (see e.g. \cite{M-M} and
\cite{Me2}) who considered symmetric
models with $AB$ interaction including the so-called Ising model
or matrices coupled in chain model, i.e $m=1,\, p=0$ and
$$P(A_1,\cdots, A_n)=\sum_{i=1}^n P(A_i) +
\sum_{i=1}^{n-1} A_{i}A_{i+1}.$$
By orthogonal polynomial techniques,
they could obtain the asymptotic behaviour of the
associated free energy when integration 
holds over Hermitian matrices.
By using completely different techniques
based on large deviations, similar 
asymptotics could be derived in \cite{GZ} and \cite{G}
for $AB$ interaction
models where the symmetry
between the matrices can be broken (i.e. we can choose $P(A_1,\cdots, A_n)=\sum_{i=1}^n P_i(A_i) +
\sum_{i=1}^{n-1} A_{i}A_{i+1}$, possibly with different $P_i$'s)
and  integration
can also hold over the orthogonal
ensemble. These techniques have moreover the advantage
to
allow the description
of the asymptotic behaviour
of the spectral measures of the matrices $(A_1,\cdots,A_n)$
with law $\mu_N^P$, key step to 
try to obtain the full expansion of $Z_N(P)$.\\
On a less rigorous
ground,
a few other models have been studied.
 The main idea to study 
 most of them 
is based on character expansion,
a technique which was introduced by A. Migdal in \cite{Mig} and 
 by C. Itzykson and J.-B. Zuber in their famous article 
on planar approximation \cite{ItZ}, 
and then widely developed in the 90's by various 
physicists (see for example \cite{DFI}, \cite{KZ} for the 
so-called $ABAB$ model or refer to \cite{K2}
for a review).
This technique allows to express the involved 
matrix integrals in terms 
basically of a sum over characters
which are simpler to deal with because the interaction is reduced
 to spherical integrals, whose asymptotics are described in \cite{GZ}.
However, this sum is in general 
an infinite signed
series (which actually might diverge),
point which is not addressed 
for instance in
  \cite{K2}. A formal 
expansion was also obtained by B. Collins 
in \cite{Co} in a very general
setting. He could obtain a formula
for the free energy of matrix
integrals as a formal series
and study the convergence 
of each terms of this series. 
However, he could not prove that
the series in fact converges.

\medskip
\def\tr{{\mbox{tr}}}

In the present article, we show
how the idea of character 
expansion can be used 
to estimate rigorously 
the   specific matrix integral in which, 
$A_N$ and $B_N$ being two $N \times N$ given Hermitian
 matrices, 
the partition function  is 
\begin{eqnarray}
 Z_N(\Phi) &\equiv &\int
 dM e^{-\frac N 2 \tr M^2 - \tr \otimes \tr \log (I \otimes I
 - B_N \otimes \Phi(M)A_N) },\label{fnpart1}\\
&=& \int
 dM e^{-\frac N 2 \tr M^2 +\sum_{k\ge 1} k^{-1} 
\tr(B_N^k) \tr((\Phi(M)A_N)^k) 
}\nonumber\\
\nonumber
\end{eqnarray}
with the following notations :
\begin{itemize}
\item $dM$ is the Lebesgue measure over 
the set $\HNC$ of Hermitian matrices of size $N$, 
\item $\tr$ is the usual trace on $\MNC$ and $I$ is the identity 
in $\MNC$, 
\item $\Phi$ is a continuous  function from $\RR$ into $\RR$.
$\Phi(M)$ is then uniquely defined by \hfill

\nn
$\Phi(M)=U\mbox{diag}(
\Phi(\l_1),\cdots,\Phi(\l_N))U^*$ when
$M=U\mbox{diag}(\l_1,\cdots,\l_N)
U^*$ for some $U\in\UNC$.
\end{itemize}
This model was studied in the case 
where $\Phi(x)=x$ in \cite{KSW}
where it was called the ``dually weighted graphs
model", because it describes, in the large $N$ limit, planar graphs
 having arbitrary coordination dependent 
weights for both vertices and faces. Note that in fact, in the case where 
$\Phi(x)=x$, the expansion is diverging (see \cite{KSW}, (2.7)).
In this work, we shall
restrict ourselves to 
functions $\Phi$ satisfying appropriate 
boundness conditions to insure 
that the  partition function $ Z_N(\Phi)$ and
its character expansion
are  well defined. 
We discuss in section \ref{KSWsec} the
relation between our result, \cite{KSW}
and the enumeration of maps.
Our main results can be sketched 
as follows

\begin{theo}
\label{maintheo}
\begin{enumerate}
\item Under appropriate assumptions   (see 
hypotheses \ref{hypothese}, 
\ref{hypothese2}),
$$F_N(\Phi)= {1\over N^2}\log Z_N(\Phi)$$
converges as $N$ goes to infinity  and a formula is derived
(see Theorem \ref{ldp}
for details).
\item Under appropriate  additional assumptions,
we can give a weak characterization
of the limit points of the
spectral measure of $M$
under the Gibbs measure associated to
$Z_N(\Phi)$ (see Proposition \ref{mini})
\end{enumerate}
\end{theo}
The main advantage
of this model
is that its character
expansion is not signed (i.e is a
sum of non negative terms),
allowing standard
Laplace method techniques.
But let us explain what we mean by ``character expansion'', i.e. expansion 
in terms of Schur polynomials.
For that, we 
recall  the following notions (see for example section 4.4. of the book \cite{Sag} for more details):
\begin{itemize}
\item a Young shape $\lambda$ is a finite sequence of non-negative integers $(\lambda_1, \lambda_2, \ldots, \lambda_l)$ written in non-increasing order. One should 
think of it as a diagram whose $i$th line is made of $\lambda_i$ empty boxes.
We denote by $|\lambda|= \sum_i \lambda_i$  the total number 
of boxes of the shape $\lambda$.\\
In the sequel, when we have a shape $\lambda = (\lambda_1, \lambda_2, \ldots)$ and an integer $N$ greater than the number of lines of $\lambda$ having a strictly positive length, we will define a sequence $l$ associated to $\lambda$ and $N$, which is an $N$-uple of integers $l_i= \lambda_i +N-i$. In particular we have that $l_1 > l_2 > \ldots > l_N \pgq 0$ and $l_i-l_{i+1}\geq 1$.

\item  for some fixed $N \in \NN$, a Young tableau will
 be any filling of the Young shape above with integers 
from $1$ to $N$ which is non-decreasing on each line and
 (strictly) increasing on each column. For each such filling, 
we define the content of a Young tableau as the $N$-uple $(\mu_1,
 \ldots, \mu_N)$ where $\mu_i$ is the number of $i$'s written in the
 tableau. \\ 
Notice that, for $N\in\NN$, a Young shape 
can be filled with integers from $1$ to $N$ if and only 
if $\l_i =0$ for $i>N$.

\item for a Young shape $\lambda$ and an integer $N$, the Schur
 polynomial $s_{\lambda}$ is an element of 
$\CC\langle x_1, \ldots, x_N\rangle$ defined by
\begin{equation}\label{defs}
 s_{\lambda}(x_1, \ldots, x_N) = \sum_T x_1^{\mu_1} \ldots 
x_N^{\mu_N},
\end{equation}
where the sum is taken over all Young tableaux $T$ of fixed 
shape $\lambda$ and $(\mu_1, \ldots, \mu_N)$ is the content of $T$.
Note that $s_\l$ is positive whenever the $x_i$'s are
and, although it is not obvious from this definition (cf for example \cite{Sag}
for a proof), $s_\l$ is a symmetric function
of the $x_i$'s.
\end{itemize}
\vspace{0.3cm}
If  $A$ is a matrix in $\MNC$, then define $s_{\lambda}(A)
 \equiv  s_{\lambda}(A_1, \ldots, A_N)$, where the 
$A_i$'s are the eigenvalues of $A$. \\
Now the point is that we shall see in Theorem \ref{formula}, whose derivation 
is the object of section \ref{cha}, that
we can write $Z_N(\Phi)$ as
$$Z_N(\Phi)=c_N\sum_\l s_\l(A_N) s_\l(B_N) Z_N(\Phi,\l)$$
where
 the sum runs over Young tableaux $\l=(\l_1\ge\l_2\cdots \ge\l_N)$
 and $ Z_N(\Phi,\l)$
is a positive function of the shape $\l$
which depends `almost continuously'
on the empirical measure $$\mun_\l:={1\over N}\sum_{i=1}^N
\d_{\l_i+N-i\over N}\in \Pa(\R^+)$$ 
where $\Pa(\R^+)$ denotes the set of probability
measures on $\R^+$.
Therefore, to study 
the asymptotic behaviour 
of $Z_N(\Phi)$ we are lead to  estimate the deviations
of more general measures and establish the following

\begin{theo}\label{scc}
Let $F:\Pa(\R^+)\ra\R$ 
be a bounded continuous 
function, and $c:\R^+\ra \R$ be a continuous
function
such that $\liminf_{ x\ra +\infty} x^{-1} c(x)>0$.
Let $(A_N, B_N)_{N\ge 0}$ be two
sequences of matrices with 
 eigenvalues taking their values 
in $[\e,1]$ for some $\e>0$ and such
that
 the spectral measures of $A_N$ and $B_N$ converge
towards $\mu_A$ and $\mu_B$  respectively.
Let $a,b\ge 0$ and
consider
the positive measure on $\Pa(\R^+)$
given, for any measurable subset 
$M\in \Pa(\R^+)$,
by 
\begin{equation}
\label{PiN}
\Pi^N(M)=\sum_{\l} 1_{\mun_\l\in M}
s_\l(A_N)^a s_\l(B_N)^b e^{N^2 F(\mun_\l)-N^2\int c(x)d\mun_\l(x)}.
\end{equation}

Then, if we   equip $\Pa(\R^+)$  with the standard weak  topology, $(\Pi^N)_{N\ge 0}$
satisfies  large deviation bounds 
with a  rate function
$H$ which is infinite on $\La^c$ where
$$\La:=\left\{\nu\in\Ppos :
d\nu(x)\ll dx,\quad {d\nu(x)\over dx}\le 1 \right\}$$
  and otherwise 
given 
by

$$H(\nu)=\int c(x) d\nu(x)-{a+b\over 2}\Sigma(\nu)
 -F(\nu) 
-aI(\log\sharp\mu_A, \nu)-b I(\log\sharp\mu_B, \nu)
-{a\over 2} S(\mu_A)-{b\over 2}S(\mu_B)$$
where 
\begin{itemize}
\item $I(\mu,\nu)$ will be defined in subsection \ref{defI},
\item $ \displaystyle \Sigma(\nu)=
\int\!\!\int \log |x-y| d\nu(x)d\nu(y) ,$\\
\item
$ \displaystyle S(\mu)=\int\!\! \int \log \left(s(x,y)
\right)d\mu(x)d\mu(y),$
with
$$s(x,y)= \int_0^1(\a x+(1-a)y)^{-1}\mbox{ if }x\neq y,\quad s(x,x)=x^{-1}\mbox{ otherwise.}$$ 
\item and for $\mu \in \PR$ and any measurable function $f : \RR \ra \RR$, we denote by 
$f \sharp \mu$ the probability measure such that, for any bounded measurable function $g$ on $\RR$,
 $ \displaystyle f \sharp \mu(g) = \int g(f(x)) d\mu(x)$.  
\end{itemize}
More precisely,

\begin{enumerate}
\item $H$ has compactly supported level sets, i.e
$\{\nu\in\Pa(\R^+) :H(\nu)\le M\}$
is compact for all $M<\infty$.
\item For any closed set $F\in\Pa(\R^+)$
$$\limsup_{N\ra\infty}
{1\over N^2}\log \Pi^N(F)\le -\inf \{ H(\nu),\nu\in F\}$$
\item For any open set $O\in\Pa(\R^+)$
$$\liminf_{N\ra\infty}
{1\over N^2}\log \Pi^N(O)\ge -\inf \{ H(\nu),\nu\in O\}$$
\end{enumerate}
In particular,

$$\lim_{N\ra\infty}
{1\over N^2}\log \Pi^N(\Pa(\R^+))= -\inf \{ H(\nu)\}$$
and the infimum is achieved.

\end{theo}
Theorem \ref{ldp}
would be a direct consequence 
of Theorem \ref{scc} (with $a=b=1$
and $\log Z_N(\Phi, \l)=N^2F(\mun_\l)-N^2\int c(x)d\mun_\l(x)$)
if $Z_N(\Phi, \l)$ was indeed a  
continuous function of $\mun_\l$ and decayed sufficiently 
fast as the size of the tableau goes to
infinity. Although it is not exactly the case, 
 most 
of the technicalities are already
contained in the proof Theorem \ref{scc},
which, as we shall see in section \ref{KSWsec},
is of independent interest.
Its proof
 relies on techniques 
developed in \cite{BA-G} in a continuous
setting, the relation
of Schur functions with spherical integrals 
(see section \ref{cha})
and  on \cite{GZ}
where the asymptotics of such integrals
were obtained.
However, the proof remains rather technical 
for various
reasons, the most severe being that
we need to define the spherical integrals
in a broader set than what was studied in \cite{GZ}.
In section \ref{proofscc}, we prove Theorem \ref{scc} 
in details.
 We precise the strategy used to show the
 Theorem \ref{scc}
 at the beginning of section \ref{proofscc}, 
just after the precise statement of the theorem.
We outline how to adapt the proofs to
obtain Theorem \ref{ldp} in section 
 \ref{lap}. 
Section \ref{secmin} is devoted to the study of the minimizers
of the rate function associated 
with the asymptotics of 
$Z_N(\Phi)$. They are reminiscent 
of \cite{KSW}
since they are described in terms
of an additional
measure describing the
optimal shape of the Young
tableau. They involve also,
following \cite{G}
and \cite{Ma},
the solutions of an Euler equation 
for isentropic flow with negative pressure
$p(\rho)=-{\pi^2\over 3}\rho^3$.\\
Finally, we comment our
result, other applications
of our techniques,  and their relations
with the problem of the enumeration of maps
in section \ref{KSWsec}.

\section{Formulation
 of the matrix model as a sum over characters}\label{cha}

Before going into the details of the large deviation principles
we have announced in the introduction, we devote this section 
to show the character expansion for $Z_N(\Phi)$ (see Theorem 
\ref{formula}). This will be useful in section \ref{lap} and can also
be seen as a justification for the definition of $\Pi^N$ we introduced above
and therefore as a motivation to prove such a result like Theorem \ref{scc}.\\  

Since we shall later also be interested by 
the Gibbs measure associated
with such a model we more generally define, after (\ref{fnpart1}),
if  $X$ is a measurable subset
of $\Pa(\RR)$ 
\begin{equation}
\label{fnpart}
 Z_N(\Phi)(X) \equiv \int_{\mun_M\in X}
 dM e^{-\frac N 2 tr M^2 - tr \otimes tr \log (I \otimes I - B_N \otimes \Phi(M)A_N) },
\end{equation}
where, for an
Hermitian matrix $M\in\Ha_N(\CC)$
with eigenvalues $(M_1,\cdots,M_N)\in\R^N$,
we shall denote $\mun_M$ the spectral measure of
$M$ given by
$$\mun_M={1\over N}\sum_{i=1}^N \d_{M_i}.$$
$\mun_M$ is an element of the space
$\Pa(\RR)$ of probability measures on the real line.
We  
endow $\Pa(\RR)$ with its usual weak topology
(i.e $\mu_n\in\Pa(\RR)$ converges towards $\mu$ iff 
$\mu_n(f)=\int f d\mu_n$ converges to $\mu(f)$
for all $f$ in the space $\Ca_b(\RR)$
of bounded continuous functions).\\
We shall assume that

\begin{hyp}\label{hypothese} 
\mbox{}
\begin{enumerate} 
\item  If $\|.\|_N$ 
denotes  the operator norm in $\MNC$, 
$\sup_{N \in \NN} \| A_N \|_N$ and  
$\sup_{N \in \NN} \| B_N \|_N$ are finite  and $\Phi$ is 
bounded. Without loss of generality,
we will assume hereafter that $$\sup_{N \in \NN} \| A_N \|_N\le 1,
\quad\sup_{N \in \NN} \| B_N \|_N\le 1$$
This amounts to multiply
$\Phi$ by $\sup_{N \in \NN} \| A_N \|_N. 
\sup_{N \in \NN} \| B_N \|_N$.
\item For all $N\in\NN$, $A_N$ and $B_N$ are non-negative
and  $\Phi$ takes its value in $\RR^+$.
\item If we define $\rho_\Phi:=-\log ||\Phi||_\infty$,
we assume that 
\begin{equation}
\label{borne}
e^{-\rho_\Phi}:=||\Phi||_\infty < 1.
\end{equation}
\end{enumerate}
\end{hyp}

Note that this assumption insures that for each $N$, 
$I \otimes I - B_N \otimes \Phi(M)A_N$ 
has positive eigenvalues, so that its logarithm is well 
defined and $\tr \otimes \tr 
\log (I \otimes I - B_N \otimes \Phi(M)A_N)$ is
bounded so that the partition function itself is well defined. 


The goal of this section is 
to express the partition function $ Z_N(\Phi)(X)$
in terms of spherical integrals, where a 
 spherical integral  $I_N$ over the unitary
group is given, 
for  two real diagonal matrices $D_N,E_N$, 
by
$$I_N(D_N,E_N):=\int \exp\{N\tr(UD_NU^*E_N)\} dm_N(U),$$
where $m_N$ denote the Haar measure  on the
 unitary group
${\cal U}_N$. 
In the sequel, we will denote $\D$ the VanderMonde determinant
given, for any diagonal
matrix $A_N=\mbox{diag}(a_1,\cdots,a_N)$, by 
$\D(A_N)=\D(a)=\prod_{i<j}|a_i-a_j|$.\\
The main result 
of this section is  

\begin{theo}\label{formula}

When Hypothesis \ref{hypothese} is satisfied, we have that

\begin{equation}
\label{Zreduc}
 Z_N(\Phi)(X)= c_N
 \sum_{\lambda} 
s_\l(A_N)s_{\lambda}(B_N)  Z_N(\Phi, \l)(X)
\end{equation}
where :

\begin{itemize}
\item $\mathcal U_N$ is the unitary group of dimension $N$,
\item the sum holds over all Young shapes, 
\item $s_\l$ is the Schur polynomial corresponding
to a Young shape $\l$, \\
\item \mbox{} \\  \vspace{-1.2cm}  
$$  \hspace{-1.3cm} Z_N(\Phi, \l)(X)= \int_{\mun_M\in X}  I_N \left(\log \Phi(M), \frac l N \right)
 \frac{\Delta(\log \Phi(M))}{\Delta(\Phi(M))}  
\Delta(M)^2 e^{-\frac N 2 \sum_{i=1}^N M_i^2} \prod_{i=1}^{N} dM_i,
$$
where $l$ is the sequence associated to $\l$ and $N$,
\item $c_N$ is a constant which only depend on $N$.
\end{itemize}
Denoting $|\l|=\sum_{i}\l_i$, we can rewrite (\ref{Zreduc})
into 
\begin{equation}
\label{Zreduc2}
 Z_N(\Phi)(X)= c_N
 \sum_{\lambda} 
s_\l(A_N)s_{\lambda}(B_N)  Z_N(\Psi, \l)(X)\, e^{-\rho_\Phi |\l|}
\end{equation}
where 
$\Psi=(||\Phi||_\infty)^{-1}\Phi$ and $c_N$ is a constant which
only depend on $N$.
\end{theo}
\prf

\begin{enumerate}
\item {\it Expansion along Young tableaux}

By definition, 
if $(B_{N,i})_{1 \ppq i \ppq N}$ 
and $((\Phi(M)A_N)_i)_{1 \ppq i \ppq N}$ are 
respectively the eigenvalues of $B_N$ and $\Phi(M)A_N$, 
we can rewrite : 
\begin{equation}
\label{explog}
e^{-tr \otimes tr \log (I \otimes I - B_N \otimes \Phi(M)A_N) } =
 \prod_{i,j=1}^N \frac{1}{1-B_{N,i}(\Phi(M)A_N)_j},
\end{equation}
where condition (\ref{borne}) ensures the
 existence of the right hand side.\\
The Cauchy formula (for a reference and a proof, 
see for example formula 4.8.4
in the book of Sagan \cite{Sag})
gives us that 
\begin{equation}
\label{Cauchy}
\prod_{i,j=1}^N \frac{1}{1-B_{N,i}(\Phi(M)A_N)_j} = 
\sum_{\lambda} s_{\lambda}(B_N) s_{\lambda}(\Phi(M)A_N),
\end{equation}
where $\lambda$ is the shape of a Young tableau and 
$ s_{\lambda}$ is the Schur polynomial corresponding 
to this shape.\\
Note that $s_\l(B_N)\ge 0$
since $B_N\ge 0$ as well as 
$s_{\lambda}(\Phi(M)A_N)=s_{\lambda}
(A_N^{1\over 2}\Phi(M)A_N^{1\over 2})\ge 0$.
Hence, the above series
converges absolutely and
we can use 
 Fubini's 
theorem to write our partition function
\begin{equation}
\label{fnpartnew}
Z_N(\Phi)(X) = \sum_{\lambda}
s_{\lambda}(B_N) \int_{\mun_M\in X} 
 e^{-\frac N 2 \tr M^2} s_{\lambda}(\Phi(M)A_N) dM.
\end{equation}

\item {\it Formulating $Z_N(\Phi)(X)$ in terms of Schur polynomials
}

It is useful to recall now the result of Weyl which establishes 
that $ s_{\lambda}$ coincides with the character 
of the unitary group associated to the shape $\lambda$ 
(this is contained in 
 theorem 7.5.B of \cite{Wey}).
 This allows us to apply to our $ s_{\lambda}$'s a key fact about
 characters : the well known property of orthogonality. 
More precisely, if $V$ and $W$ are two unitary matrices of 
size $N$, this property reads, for any shape $\lambda$,
\begin{equation}
\label{ortho1}
\int  s_{\lambda}(UVU^*W) dm_N(U)= 
\frac{1}{d_{\lambda}} s_{\lambda}(V)s_{\lambda}(W),
\end{equation}
where $dm_N$ is the Haar 
measure on the unitary group $\mathcal U_N$ normalized to have mass 
one and $d_{\lambda}=s_{\lambda}(1,1,\cdots ,1)$. Its explicit form is 
\begin{equation}
\label{diml}
d_{\lambda} = 
\frac{\Delta(l)}{\prod_{i=1}^{N-1} i!},
\end{equation}
with $l = \mbox{diag}(l_1, \ldots, l_N)$ where we recall that $l_i = \l_i +N-i$.\\
A proof of formula (\ref{ortho1}) can be easily deduced from proposition $\textrm{I\!I}$.4.2 of \cite{BtD} (see also exercise 3 p.84 therein) whereas the explicit expression of $d_{\lambda}$ given in (\ref{diml}) appears in \cite{Wey}.\\

As a consequence, with the notations 
introduced above,
\begin{equation}
\label{ortho}
\int  s_{\lambda}(U\Phi(M)U^*A_N) dm_N(U)= 
\frac{1}{d_{\lambda}} s_{\lambda}(\Phi(M))s_{\lambda}(A_N).
\end{equation}
Combining equations (\ref{fnpartnew}) and 
(\ref{ortho}), we can rewrite our partition function
\begin{equation}
\label{Zsimplif}
Z_N(\Phi)(X) = c^\prime_N \sum_{\lambda} \frac{1}{d_{\lambda}} 
s_{\lambda}(A_N)s_{\lambda}(B_N) 
\int_{\mun_M\in X} s_{\lambda}(\Phi(M)) e^{-\frac N 2 \tr M^2}
 \Delta(M)^2 \prod_{i=1}^N dM_i,
\end{equation}
where $\prod_{i=1}^{N} dM_i$ is the product Lebesgue 
measure on $\RR^N$ and $c^\prime_N$ some normalizing 
constant, only depending on $N$. \\

\item{\it Relation between Schur polynomials and spherical integrals}

We can now recall the following determinantal formula 
for $ s_{\lambda}$, that can be found for example 
in corollary 4.6.2 of  \cite{Sag}:
\begin{equation}
\label{Weyl}
s_{\lambda} (\bf x \rm ) = 
\frac{det(x_i^{l_j})_{i,j}}{\Delta(\bf x \rm)},
\end{equation} 
where $\Delta$ is the VanderMonde determinant, 
$\mathbf{x} = (x_i)_{1 \ppq i \ppq N}$ and $l$ is the tableau associated to $\lambda$ (that is to say $l_j= \l_j+N-j$ for $1 \le j \le N$). \\


We then use a formula due to Harish-Chandra (see \cite{Meh}): 
if $C_N$ and $D_N$ are two $N \times N$ matrices whose eigenvalues 
$C_N(i)$ and $D_N(j)$ are distinct, we have that
\begin{equation}
\label{IZ}
I_N(C_N, D_N) = \frac{det(\exp NC_N(i)D_N(j))_{i,j}}{\Delta(C_N)\Delta(D_N)}.
\end{equation}


This last equation together with the determinantal formula (\ref{Weyl}) allows
 us to rewrite
for any $M\in\HNC$ with non negative distinct eigenvalues: 

\begin{equation}
\label{IZHC}
s_{\lambda}(M) = I_N \left(\log M, 
\frac l N  \right) \Delta \left(\frac l N \right)
 \frac{\Delta(\log M)}{\Delta(M)},
\end{equation}
Note that under the measure $ e^{-\frac N 2 \tr M^2} dM$, 
the eigenvalues of the matrix $M$ 
are almost surely distinct, and therefore 
so are the eigenvalues 
of the two matrices $\Phi(M)$ and $\log \Phi(M)$ 
by hypothesis \ref{hypothese}.3. Note however that
(\ref{IZHC}) extends readily to any non negative matrix by extending the definition
$$\frac{\Delta(\log M)}{\Delta(M)}=e^{\sum_{i<j}s(\l_i,\l_j)},$$\\
with $s$ as defined in Theorem \ref{scc}.\\
From (\ref{IZHC}), we conclude 
that there exists a constant $c_N$ depending only on $N$ such that,  
\begin{multline*}
Z_N(\Phi)(X) = c_N
\sum_{\lambda} s_{\lambda}(A_N)s_{\lambda}(B_N)\\
\ts
\int_{\mun_M\in X}  I_N \left(\log \Phi(M), \frac l N \right)
 \frac{\Delta(\log \Phi(M))}{\Delta(\Phi(M))}  
\Delta(M)^2 e^{-\frac N 2 \sum_{i=1}^N M_i^2} \prod_{i=1}^{N} dM_i, 
\end{multline*}
which completes the proof of Theorem \ref{formula}
except from formula (\ref{Zreduc2})
which is easily obtained by dividing the $\Phi$
by its norm before beginning the expansion.
\end{enumerate}

\hfill\xx

\begin{rmk}
\label{cN}
If we denote by
$$\mbox{vol}(\mathcal U_N) := \frac{\int  e^{-\frac N 2 \tr M^2}
  dM}{\int  e^{-\frac N 2 \tr M^2}
 \Delta(M)^2 \prod_{i=1}^N dM_i},$$
 we can easily deduce
from equations (\ref{diml}), (\ref{Zsimplif}) and (\ref{IZHC})
above, that our normalizing constant $c_N$ is given by
$$ c_N = \mbox{vol}(\mathcal U_N)\left( {\prod_{i=1}^{N-1} i! \over N^{N(N-1) \over 2}}
 \right). $$
\end{rmk}

\section{Large deviations estimates for the empirical
distribution \\
of Young tableaux following the law $\Pi^N$}
\label{proofscc}
The object of this section is to prove Theorem \ref{scc}.\\
From the definition (\ref{PiN}) and following (\ref{IZHC}),
we get that $\Pi^N$ is the positive measure given,
for any measurable subset $M$ of $\Ppos$, by :
\begin{multline*}
\Pi^N(M)=
e^{{a\over 2}N^2S_N(\mun_A)+{b\over 2}N^2S_N(\mun_B)}\\ 
\ts \sum_{\l:\mun_\l\in M} \D\left({l\over N}\right)^{a+b}
I_N\left(\log A_N, {l\over N}\right)^aI_N\left(\log B_N, {l\over N}\right)^b
 e^{N^2 F(\mun_\l)-N^2\int c(x)d\mun_\l(x)}
\end{multline*}
where
$$e^{{N^2\over 2}S_N(\mun_A)}:={\D(\log(A_N))\over \D(A_N)}.$$
Let us denote
$$\tilde\Pi^N(M)=
\sum_{\l} 1_{\mun_\l\in M} \D\left({l\over N}\right)^{a+b}
I_N\left(\log A_N, {l\over N}\right)^aI_N\left(\log B_N, {l\over N}\right)^b
 e^{N^2 F(\mun_\l)-N^2\int c(x)d\mun_\l(x)}
.$$
We shall prove in this section

\begin{theo}\label{scc2}
Let $(F,c, (A_N,B_N),a,b)$
be as in Theorem \ref{scc}.\\
Then $(\tilde\Pi^N)_{N\ge 0}$
satisfies  large deviation bounds 
with  rate function
$\tilde H$ 
which is infinite on $\La^c$
 and otherwise 
given 
by

$$\tilde H(\nu)=\int c(x)d\nu(x)-{a+b\over 2}\Sigma(\nu)
 -F(\nu)
-aI(\log\sharp\mu_A, \nu)-b I(\log\sharp\mu_B, \nu).
$$
More precisely,

\begin{enumerate}
\item 
$\{\nu\in\Pa(\R^+) : \tilde H(\nu)\le M\}$
is compact for all $M<\infty$.
\item For any closed set $F\in\Pa(\R^+)$,
$$\limsup_{N\ra\infty}
{1\over N^2}\log \tilde \Pi^N(F)\le -\inf \{ \tilde H(\nu),\nu\in F\}$$
\item For any open set $O\in\Pa(\R^+)$,
$$\liminf_{N\ra\infty}
{1\over N^2}\log \tilde \Pi^N(O)\ge  -\inf \{\tilde H(\nu),\nu\in O\}$$
\end{enumerate}

\end{theo}
Theorem
\ref{scc} is easily deduced from
Theorem \ref{scc2} since

\begin{equation}
\label{fnpartlbw}
S_N(\mun_A)
={2\over N^2}\sum_{i< j}s(A_i,A_j).
\end{equation}
Hence,
since 
$s
$
is a  bounded continuous
function on $[\e,1]^2$,
we deduce (see Lemma 7.3.12 in \cite{DZ}) that, as $\mun_A$ converges
to $\mu_A$, 
$$\lim_{N\ra\infty}S_N(\mun_A)
=S(\mu_A)$$
and similarly for $B_N$.

\medskip

The proof of Theorem \ref{scc2} 
is heuristically simple since it amounts
to perform a  Laplace
method and notice
that 
the uniform measure on Young shape
will not produce any entropy
on the scale $N^2$. On a rigorous
ground, it becomes a bit technical, 
for mainly the two following reasons  :

\begin{itemize} 
\item The law of $\mun_\l$ is discrete
so that the arguments developed in \cite{BA-G}
to obtain large deviation principles in similar
scales and potentials
have to be adapted. In particular, the discrete nature
of the Young tableaux implies that $\tilde H$ is infinite
on $\La^c$.
\item More cumbersome is the fact that the natural space 
where the empirical measure of the Young
tableaux lives is $\Pa_1(\RR^+):=\{\nu\in\Ppos : \int x d\nu(x)<\infty\}$.
Hence, all the limiting spherical integrals appearing 
are of the type $I(\mu,\nu)$ with
$\mu$ in the set $ \Pa_\infty(\RR)$
of compactly supported probability measures
but $\nu\in \Pa_1(\RR^+)$.  Such limits were
not proved to exist
in \cite{GZ} (where $\nu(x^2)<\infty$
was assumed), the formula obtained in \cite{GZ}
is not valid,  and  continuity
statements for $I$ are lacking a priori.
 
\end{itemize}

The proof
  nevertheless
follows the usual scheme :

\begin{enumerate}
\item In subsection \ref{grf} we
study the rate function and prove that its level sets are compact.
\item In subsection \ref{exptight} we
show that the family of measures $(\tilde\Pi^N)_{N \in \NN}$ 
is exponentially tight.
More precisely, if we let $\Ka_L$ be the compact subset

$$\Ka_L=\left\{ \nu\in \Ppos :\int xd\nu(x)\le 
L\right\}$$
we prove that
$$\limsup_{L\ra\infty}\limsup_{N\ra\infty}
{1\over N^2}\log \tilde \Pi^N(\Ka_L^c)=-\infty.$$
\item In subsection \ref{ldub} we prove 
the upper bound for arbitrarily small balls, i.e
if $d$ is a metric on $\Pa(\RR)$ 
compatible with the weak topology
such as
the Dudley's metric $d$ given by
$$ d(\mu, \nu) = \sup \left|\int f d\mu - \int f d\nu \right|, $$
where the supremum is taken over all Lipschitz functions $f$
 with Lipschitz norm less than $1$ (note that this distance is compatible
 with the weak topology), and if we set
$$B(\nu,\d)=\{\mu\in\Ppos ;d(\mu,\nu)<\d\}$$
we show that for any $\nu
\in \cup_{L\in\NN} \Ka_L$,
$$\limsup_{\d\ra\infty}
\limsup_{N\ra\infty}
{1\over N^2}\log \tilde\Pi^N( B(\nu,\d))
\le -\tilde H (\nu).$$
\item In subsection \ref{ldlb} we prove 
the lower bound for arbitrarily small balls,
i.e that for any $\nu\in \cup_{L\in\NN} \Ka_L$,
$$\liminf_{\d\ra\infty}
\liminf_{N\ra\infty}
{1\over N^2}\log \tilde\Pi^N( B(\nu,\d))
\ge -\tilde H (\nu).$$
\end{enumerate}
By Theorem 4.1.11 in \cite{DZ}, the above results
prove Theorem \ref{scc2}.

\subsection{ $\tilde H $ has compact level sets}\label{grf}

To prove that $\tilde H $  has compact level sets,
we shall first define it properly, 
that is define appropriately the limit of the spherical integrals.

\subsubsection{Definition and properties of I}\label{contsph}
 Let us remind that
it was proved in 
theorem 1.1 of 
\cite{GZ} that
\begin{equation}
\label{def:I}
I(\mu_D, \mu_E):=\lim_{N \rightarrow \infty} \frac{1}{N^2} 
\log I_N(D_N, E_N)  
\end{equation}
exists for all  sequences 
of diagonal matrices $(D_N,E_N)_{N\in\NN}$
with spectral measures converging towards $\mu_D$ and $\mu_E$
respectively and such that $\sup_N||D_N||_N$ and 
 $\sup_N\mun_E(x^2)$ are finite. A formula for
$I$ is given in \cite{GZ}
when either $\Sigma(\mu_E)$ or
$\Sigma(\mu_D)$ are finite. If they are not,
the limit still exists since  spherical integrals 
are uniformly continuous (see Lemma \ref{defI}.4))
and the measures with finite $\Sigma$ are dense, but
its formula is far from being clear (see a discussion in \cite{GZ2}).
However, let us remark that  since the spherical integrals
under considerations are always bounded, the
rate function $\tilde H (\nu)$ is infinite
unless $\nu$ has finite entropy
$\Sigma$ (see the end of section \ref{grf})
so that we can always use the formula given in \cite{GZ}.


Since $\tilde H(\nu)$ is infinite if $\int x d\nu(x)=+\infty$ (see section \ref{grf})
and $\mu_A$ and $\mu_B$ are supposed to be supported on $[\e, 1]$,
 it is enough to extend the 
definition of  $I(\mu,\nu)$ to  compactly supported measures
$\mu$ with support in $\RR^-$ 
but $\nu\in \Pa_1(\RR^+)$.
We shall
prove 


\begin{lem}\label{defI}
Let $R\in\RR^+$ 
and  $\mu$ be a probability measure on $[-R,0]$
and $\nu\in \Pa_1(\RR^+)$. Then
\begin{enumerate}
\item Let $\phi_M(x)=x\wedge M$. $I(\mu,\phi_M\sharp \nu)$
is well defined and decreases towards a limit 
$$I(\mu,\nu):=\lim_{M\ra\infty} I(\mu,\phi_M\sharp \nu).$$
Moreover,  for any $M\ge 0$,
$$  I(\mu,\phi_M\sharp \nu)-R \nu( x-\phi_M(x))
\le I(\mu,\nu)\le I(\mu,\phi_M\sharp \nu).$$
\item Let $\Pa^R(\RR)=\{\mu\in\Pa(\RR):\mu([-R,R]^c)=0\}$
and $\Pa_q(\RR^+)=\{\mu\in\Pa(\RR^+):\mu(|x|^q)\le R\}$.
Then there exists a function $\kappa(\d,R)$
such that for any $R<\infty$, $\kappa(\d,R)$
goes to zero
as $\d$ goes to zero and for any
$(\mu,\mu')\in \Pa^R(\RR)$
any $(\nu,\nu')\in \Pa_2(\RR)$, such that $d(\mu, \mu') + d(\nu, \nu') < \d$,
$$ |I(\mu,\nu)-I(\mu',\nu')|\le \kappa(\d,R).$$
\item For any $\mu\in\Pa(\RR^-)$ and $\nu\in \Pa_1(\RR^+)$,
$$\mu(x)\nu(x)\le I(\mu,\nu)\le 0.$$
\item For any sequence $(D_N,E_N)$ 
of diagonal Hermitian matrices with $D_N\le 0$ and $E_N\ge 0$,
for any $M\in\RR^+$,

\begin{equation}\label{app1}
I_N(D_N,\phi_M(E_N))e^{-N||D_N||_N\tr(E_N-\phi_M(E_N))}
\le I_N(D_N,E_N)\le I_N(D_N,\phi_M(E_N)).
\end{equation}
Moreover there exists a function $g:[0,1]\ts \RR^+ \mapsto \RR^+$,
 depending on the limiting measures $\mu_E,\, \mu_D$ 
only, such that $g(\d,M)$ goes to zero
as $\d$ does for any $M\in\R^+$,
and  so that

\begin{equation}\label{cont1}
 \left|\frac{1}{N^2} \log 
\frac{I_N(\hat D_N, \phi_M(\hat E_N))}{I_N
( D_N,\phi_M(E_N))}
\right|\leq g(\delta,M)\,.
\end{equation}
for any $N\in\NN$ and any diagonal matrices $(D_N,E_N, \hat D_N,
\hat E_N)$ such that $E_N,\hat E_N$
are non-negative and 

$$ d(\mun_{D_N},\mun_{\hat D_N})+d(\mun_{E_N},\mun_{\hat
E_N})<\d,
\qquad \mun_{E_N}(x^2)+\mun_{\hat E_N}(x^2)\le M.$$
\end{enumerate}
\end{lem}
\prf

$\bullet$
We first prove the last point.
If we denote $D_N=\mbox{diag}(d_1,\cdots,d_N)$
and $E_N=\mbox{diag}(e_1,\cdots,e_N)$,
\begin{eqnarray*}
I_N(D_N,E_N)&=&  \int e^{N\tr(D_NUE_NU^*)} dm_N(U)\\
&=& \int e^{N\sum_{i,j=1}^N d_i e_j |u_{ij}|^2} dm_N(U)\\
&\le& \int e^{N\sum_{i,j=1}^N d_i \phi_M(e_j) |u_{ij}|^2} dm_N(U)\\
\end{eqnarray*}
where we used that $d_i\le 0$. 
The opposite inequality of (\ref{app1}) is also trivial 
since
\begin{eqnarray*}
I_N(D_N,E_N)&\ge& e^{N||D_N||_N\sum_{i,j=1}^N
(e_j-\phi_M(e_j))}
\int e^{N\sum_{i,j=1}^N d_i \phi_M(e_j)  |u_{ij}|^2} dm_N(U)\\
&=& e^{-N||D_N||_N \tr(E_N-\phi_M(E_N))}I_N(D_N,\phi_M(E_N))
\end{eqnarray*}
The continuity statement (\ref{cont1})  is a direct consequence 
of Lemma 5.1 in \cite{GZ}  since $\phi_M(E_N)$
is uniformly bounded by $M$ and $d(\phi_M\sharp \mu,\phi_M\sharp \mu')
\le d(\mu,\mu')$ for any $\mu,\mu'\in\Pa(\RR)$.

\medskip

$\bullet$ We can now prove the first point.
From (\ref{app1}), we deduce
that for any $M\in\RR^+$, any $E_N\ge 0$
with spectral measure converging towards
$\mu_E$ and any sequence of bounded  non-positive 
diagonal
matrices $D_N$ with spectral measure 
converging towards $\mu_D$
\begin{eqnarray}
\limsup_{N\ra\infty}
{1\over N^2}\log I_N(D_N,E_N)&\le&
\limsup_{N\ra\infty} {1\over N^2}\log I_N(D_N,\phi_M(E_N)) \nonumber\\
&=& I(\mu_D,\phi_M\sharp\mu_E),\label{ineq1}
\end{eqnarray}
where the last equality comes from
the observation that $(\phi_M(D_N),E_N)$ 
are uniformly bounded by hypothesis so that the
convergence holds by theorem 1.1 in \cite{GZ}. 
With $\mu_E=\phi_{L}\sharp\nu$ for some $L\ge M$
and $E_N$ chosen so that $\mun_{E_N}(|x|>L)=0$, 
the left hand side of 
(\ref{ineq1}) converges towards $I(\mu_D,\phi_L\sharp\nu)$
showing that 
 $M\ra I(\mu_D,\phi_M\sharp\mu_E)$ is non-increasing.
Hence, it converges towards some limit (maybe infinite at
this stage).
Now, we choose a special sequence 
$(E_N)_{N\in\NN}$ such that
$$\lim_{N\ra \infty} {1\over N}\tr(E_N-\phi_M(E_N)) =\mu_E(x-\phi_M(x)).$$
We can construct it as follows ; 
assume first that $\mu_E$ has no atoms
and set
\begin{eqnarray*}
E_{1,N} & = & \inf \left\{ x \,\, / \,\, 
\mu_E((-\infty, x]) \pgq \frac{1}{N+1}\right\} \\
E_{i+1,N} & = & 
\inf \left\{ x \pgq E_{i,N} \,\, / \,\, \mu_E((E_{i,N}, x])
 \pgq \frac{1}{N+1}\right\}.
\end{eqnarray*}
Then it is not hard to see that $\mun_{E_N}={1\over N}\sum_{i=1}^N 
\d_{E_{i,N}}$ converges towards $\mu_E$.
Moreover, 
\begin{eqnarray*}
\mun_{E_N}(x-\phi_M(x))&=&{1\over N}\sum_{E_{i,N}\ge M}
(E_{i,N}- M)\\
&\le& {N+1\over N}\sum_{E_{i,N}\ge M}
(E_{i,N}- M) \mu_E([E_{i,N},E_{i+1,N}])
\le {N+1\over N} \mu_E( (x-M)1_{x\ge M}).\\
\end{eqnarray*}
If $\mu_E$ has atoms, we consider
a finite collection of atoms $\{a_1,\cdots,a_K\}$
 such that each of the remaining atoms 
has mass smaller than $(N+1)^{-1}$.
 Then, $E_N$ has $\lfloor N\mu_E(\{a_i\}) \rfloor$
eigenvalues equal to $a_i$ for $1\le i\le K$.
The remaining eigenvalues are chosen as above.

Inequality (\ref{app1})
yields with this choice 

$${1\over N^2}\log I_N(D_N,E_N)\ge I_N(D_N,\phi_M(E_N)
) e^{-N(N+1)\sup_N ||D_N||_N \mu_E( (x-M)1_{x\ge M})}$$
and therefore
\begin{eqnarray}
\liminf_{N\ra\infty}
{1\over N^2}\log I_N(D_N,E_N)&\ge& 
-\sup_N ||D_N||_N \mu_E( (x-M)1_{x\ge M})+I(\mu_D,\phi_M\sharp\mu_E)
\label{ineq2}\\
\nonumber
\end{eqnarray}
(\ref{ineq1}) and (\ref{ineq2})
shows that for such a sequence
\begin{eqnarray}
-\sup_N ||D_N||_N \mu_E( (x-M)1_{x\ge M})+I(\mu_D,\phi_M\sharp\mu_E)
&\le& \liminf_{N\ra\infty}
{1\over N^2}\log I_N(D_N,E_N)\nonumber\\
&\le&\limsup_{N\ra\infty} {1\over N^2}\log I_N(D_N,E_N)\nonumber\\
&\le& I(\mu_D,\mu_E)\label{boulduc}\\
\nonumber
\end{eqnarray}
This completes the proof of the first point.

$\bullet$ The second point is a direct consequence 
of the fourth too.
 Indeed, let $(\mu,\mu',\nu,\nu')$ be
such that
$$d(\mu,\mu')+d(\nu,\nu')<\d.$$
Then, we  choose 
a sequence $(D_N,E_N)$  (resp. $(\hat D_N,\hat E_N)$)
of matrices 
with spectral measure converging towards 
$(\mu,\nu)$ (resp. $(\mu',\nu')$) such that
$$\max\{d(\mun_{D_N},\mu),d(\mun_{\hat D_N},\mu'),d(\mun_{ E_N},\nu)
, d(\mun_{\hat E_N},\nu')\} <{\d}$$
which implies 
$$d(\mun_{D_N},\mun_{\hat D_N})<2\d,\quad d(\mun_{ E_N}, \mun_{\hat
E_N})<2\d$$
so that 4. implies, by taking the limit as $N$ goes to infinity
(here $M=R$),
that 
$$|I(\mu,\nu)-I(\mu',\nu')|\le g(2\d, R).$$

$\bullet$ In point 3., the upper bound on $I$ is trivial
and the lower bound comes from Jensen's inequality
which yields

\begin{eqnarray*}
I_N(D_N,E_N)&=&\int e^{N\sum_{i,j=1}^N e_id_j |u_{ij}|^2} dm_N(U)\\
&\ge& e^{N\sum_{i,j=1}^N e_id_j\int  |u_{ij}|^2 dm_N(U)  }\\
&=&e^{\sum_{i,j=1}^N e_id_j}=e^{N^2\mun_{E_N}(x)\mun_{D_N}(x)}\\
\end{eqnarray*}
The result is then obtained by letting $N$ going to infinity.

\hfill\xx

\subsubsection{$\tilde H $ has compact level sets }\label{IAgrf}
In this section, we prove Theorem \ref{scc2}.1
by proving first
that $\tilde H $ is lower semi-continuous 
and then that its level sets are compact.

$\bullet$ $\tilde H $ is lower semi-continuous, i.e
$\{\nu\in\Ppos:
\tilde H (\nu)\le M\}$ is closed for any $M\in\R^+$.
We recall that $\La$ is the set of probability measures which
are absolutely continuous with respect
to Lebesgue measure and with density bounded
by one and 
note that $\{\tilde H \le M\}=\La\cap\{\tilde H ^1\le M\}$
where   $\tilde H ^1(\nu)$ is given by the same 
formula than  $\tilde H (\nu)$
even for $\nu\in \La^c$.
We first check that $\La$
is closed and then show
that $\tilde H ^1$  is lower semi-continuous,
these two points proving that $\{\tilde H \le M\}$ is closed.

\nn
To show that
$\La$ is closed,
take a sequence $(\nu_n)_{n \in \NN}$  of measures in $\La$ 
converging weakly to a measure $\nu$. For any $c$ and $d$, the 
function $\mathbf{1}_{[c,d]}$
is upper semi-continuous so that  
$$ |d-c| \ge \limsup_{n \ra \infty} \nu_n([c,d]) \ge \nu([c,d]).$$
so that $\nu$ is in $\La$.\\

We now show that
$\tilde H^1$ is a supremum of continuous
functions which we define as follows:
we let, with $\phi_M(x)=x\wedge M$ for $M\ge 0$
as in Lemma \ref{defI}, and
 for $\nu\in \Ppos$, 

$$\tilde H^M(\nu) := 
 -aI(\log \sharp \mu_A,\phi_M\sharp
 \nu) -b I(\log \sharp \mu_B,\phi_M\sharp \nu)
 +\int\!\!\int
 g(x,y)
\wedge M d\nu(x)d\nu(y) -F(\nu) $$
with
\begin{equation}\label{defg}
g(x,y)=\left({a+b\over 2}\right)\log|x-y|^{-1} +{1\over 2} c(x)+{1\over 2} c(
y)\end{equation}
We claim that for any finite $M$, $\tilde H^M$
is   continuous on $\Ppos$. Indeed,
by Lemma \ref{defI}.2,
 for $C=A$ or $B$,
 $\nu\in\Ppos\mapsto I(\log \sharp \mu_C,\phi_M\sharp
 \nu)\in\RR$  is continuous since $\log \sharp \mu_C$ is compactly
supported by hypothesis \ref{hypothese}.1.
Moreover, it is not hard to check that
$g$ is bounded below and
continuous except when on the diagonal $\{x=y\}$
where it goes to infinity. Consequently,
$g\wedge M$ is a bounded 
continuous function on $\RR^2$.
Thus
$\mu\ra \int\!\!\int g(x,y)\wedge M d\mu(x)d\mu(x)$
is  bounded continuous.

This last argument finishes to prove that
$\tilde H ^M$ is a  continuous
function on $\Ppos$. To 
deduce that $\tilde H^1 $ is lower 
semi-continuous, it is therefore enough 
to prove that

\begin{equation}\label{toto}
\tilde H^1 (\nu)=\sup_{M\ge 0} \{ \tilde H ^M(\nu)\}.
\end{equation}
But this is straightforward since
monotone convergence theorem asserts
that for any $f$ bounded below 
$$\lim_{M\uparrow \infty}\int\!\!\int f(x,y)\wedge M
d\mu(x)d\mu(y)= \int\!\!\int f(x,y)
d\mu(x)d\mu(y)$$
and by Lemma \ref{defI}.1,
$I(\mu,\phi_M\sharp\nu)$ decreases towards
its limit $I(\mu,\nu)$.

\nn
$\bullet$ As a consequence of the last point,
for any $M\ge 0$, $\{\nu\in \Ppos :
\tilde H (\nu)\le M\}$ is closed.
We now check that it is compact by 
showing that it is contained 
in a compact set. In fact, by Lemma \ref{defI}.3,
\begin{equation}\label{bb}
\tilde H (\nu)\ge \int \!\!\int g(x,y) d\nu(x)d\nu(y) -\sup_{\nu\in\Ppos} F(\nu)
\end{equation}
and it is not hard to check that,
 since
we assumed  $\liminf x^{-1}c(x)>0$,  there exists 
a finite constant $C$ and $\rho>0$
such that for any $(x,y)\in (\RR^+)^2$
\begin{equation}\label{borg}
g(x,y)\ge \frac \rho 2 x+
\frac \rho 2 y+C 
\end{equation}
yielding with (\ref{bb})
that for any $M\in\RR^+$, if $C'=C-\sup_{\nu\in\Ppos} F(\nu)$,
$$\{\nu\in \Ppos :\tilde H (\nu)\le
M\} \subset
\left\{\nu\in\Ppos : \int x d\nu(x)\le {2\over\rho} (M-C')
\right\}:=\Ka_{M,\rho}.$$
Since $\Ka_{M,\rho}$ is a compact subset
of $\Ppos$, the proof is completed.\\
Note that since $\int \!\!\int g(x,y)d\nu(x)d\nu(y)=\int c(x)d\nu(x)
-\Sigma(\nu)$ and $c$ is bounded below,  we also see from (\ref{bb})
that $\tilde H (\nu)<\infty$ implies
 $|\Sigma(\nu)|<\infty$.

\hfill\xx

\subsection{$\tilde\Pi^N$ is exponentially tight}\label{exptight}
The goal of this section
is to prove that
\begin{lem}\label{exptightlem}
$\tilde \Pi^N$ is exponentially tight, and more precisely if
we set
$$\Ka_L:= \left\{\nu\in\Ppos : \int x d\nu(x)\le 
L
\right\},$$
then
$$\limsup_{L\ra\infty}\limsup_{N\ra\infty}
{1\over N^2}\log
\tilde\Pi^N(\Ka_L^c)=-\infty.$$
\end{lem}
\prf \ 
Since the spherical integrals under consideration
are uniformly bounded above by one and $F$ is uniformly bounded
by a constant $||F||_\infty$,
$$\tilde \Pi^N(X)\le e^{N^2 ||F||_\infty}
 \sum_{\lambda :\mun_\l\in X } 
 e^{-{N^2}\int_{x\neq y}g(x,y)d\mun_\l(x)d\mun_\l(y)
},$$
Choosing $X=\Ka_L^c$, we get by (\ref{borg})
that
\begin{equation}\label{jea}
\tilde \Pi^N(\Ka_L^c )\le e^{N^2 ||F||_\infty +N^2 C}
 \sum_\l 1_{ \mun_\l\in \Ka_L^c}
e^{-N^2\rho\int x d\mun_\l(x)}\D\left({l\over N}\right)^{a+b}\end{equation}
It remains to consider the sums over Young shapes.
Let us recall that 
$$\mun_\l(x)={1\over N^2}\sum_{i=1}^N l_i
={1\over N}\sum_{i=1}^N \left({\l_i\over N}-{i\over N}\right)+1 \le 
N^{-2}|\l|_N   + 1 $$
where $|\l|_N=\sum_{i\le N}\l_i$.
Therefore, for any $L\ge 0$,
\begin{eqnarray*}
\sum_{\lambda :\mun_\l(x)\ge L } 
 e^{-\rho|\l|}\D\left({l\over N}\right)^{a+b}
 & \le & \sum_{\lambda :|\l|_N \ge N^2( L-1)}e^{-\rho |\l|}
\D\left({l\over N}\right)^{a+b} \\
& \le & e^{- \frac 1 2 \rho N^2(L-1)}
 \sum_{\lambda :|\l|_N \ge N^2( L-1)}
e^{-\frac 1 2 \rho|\l|}\D\left({l\over N}\right)^{a+b}.
\end{eqnarray*}
For any $j$,
$$ \prod_{j<i} \left|{l_i\over N} - {l_j\over N}\right|
 \le \left( {l_j\over N}\right)^{N-j},$$
therefore, for any shape,
$$\D\left({l\over N}\right)^{a+b} 
e^{- \frac 1 2 \rho |\l|} \le e^{(a+b)\sum_j (N-j)
 \log {l_j\over N} - \frac 1 4 N \rho {l_j\over N}} \le e^{N^2 C^{''}},$$
where $ C^{''} = \sup_x\left\{(a+b)\log x - \frac 1 4  \rho x\right\}- \frac 1 8$.\\
Now the number of Young shapes $\l$
such that $|\l|_N= m$ is bounded by $ C_m^N$ so that we
conclude
\begin{eqnarray}
\sum_{\lambda :|\l|_N \ge N^2( L-1)} \hspace{-0.8cm}
e^{-\rho |\l|}\D\left({l\over N}\right)^{a+b}
&\le& e^{N^2 C}e^{- \frac 1 2 \rho N^2(L-1)}{1\over N!}
\sum_{m\ge N^2( L-1)}  \hspace{-0.5cm} m(m-1)\cdots (m-N+1) 
e^{-\frac 1 4 \rho m},\nonumber\\
&\le& e^{N^2 C}e^{- \frac 1 2 \rho
 N^2(L-1)} {1\over N!} \sum_{m\ge N^2(L-1)}
 \hspace{-0.5cm}e^{N\log m}e^{-\rho m}\nonumber\\
&\le & e^{- \frac 1 2 (\rho-\d) N^2(L-1)} \label{b2}\\
\nonumber
\end{eqnarray}
where in the last line $\d$
is any positive number and
the inequality holds as soon as
$N$ and $L$ are big enough.
 (\ref{jea})
and (\ref{b2}) give Lemma \ref{exptightlem}.

\hfill\xx

\subsection{$(\tilde\Pi^N)_{N\ge 0}$ satisfies a weak large deviation upper
bound}\label{ldub}
In this section, we shall prove the following

\begin{lem}
 $\tilde\Pi^N$ satisfies a weak
large deviation upper bound 
in the scale $N^2$  with  rate function
 $\tilde H$
i.e for any $\nu\in \Ppos$,
$$\limsup_{\d\ra 0}
\limsup_{N\ra\infty} {1\over N^2}\log \tilde\Pi^N(
B(\nu,\d))
\le - \tilde H(\nu)
.$$
\end{lem}
\prf \
We first prove that for any $\eps >0$, if $\nu$ is such that there exists two positive real
numbers $\a$ and $\b$ $(\a<\b)$ such that
$\nu([\a,\b]) \ge (1+\eps)(\b-\a)$, then,
\begin{equation}
\label{horsL}
 \limsup_{\d\ra 0}
\limsup_{N\ra\infty} {1\over N^2}\log \tilde\Pi^N(
B(\nu,\d))
= - \infty.
\end{equation}
The main remark is that, for any shape $\l$, as the $l_i$ are (strictly) decreasing we have that, for any $c<d$,
\begin{eqnarray}
\label{mudc}
\mun_\l([c,d]) & = & \frac 1 N \sharp\left\{ i : \frac{l_i}{N} \in [c,d] \right\} 
                \le  \frac 1 N (\lfloor N(d-c)\rfloor +1) \nonumber\\
		   & \le & \left(1 + \frac \eps 2 \right)(d-c),
\end{eqnarray}
where the last inequality holds for $N$ large enough.\\
Let be $\eta >0$ and consider the function $f: \RR \ra \RR$ such that
\[
f(x) = \left\{\begin{array}{ll}
0, & \mbox{if $x<\a-\eta$ or $x>\b+\eta$,}\\
\frac 1 2 (x-\a-\eta),& \mbox{if $\a-\eta \le x \le \a$,}\\
\eta ,& \mbox{if $\a<x<\b$,}\\
\frac 1 2 (-x+\b+\eta),& \mbox{if $\b \le x \le \b+ \eta$.}\\

\end{array}
\right.
\]
Note that, for $\eta$ small enough, the Lipschitz
 norm of $f$ is bounded by $1$.\\
And we have, for any shape $\l$,
\begin{eqnarray*}
\int f d\nu - \int f d\mun_\l & = & \int_{\a-\eta}^{\a} f (d\nu -  d\mun_\l ) +\int_{\b}^{\b+\eta} f( d\nu -  d\mun_\l ) + \int_{\a}^{\b} f (d\nu -d\mun_\l ).
\end{eqnarray*}
Using (\ref{mudc}) twice, we get that, for any shape $\l$ and $N$ large enough,
$$ \int_{\a-\eta}^{\a} f d\mun_\l \le \frac{\eta^2}{2},$$
(and the same thing for $\b$) and that
$$ \int_{\a}^{\b} f d\nu - \int_{\a}^{\b} f d\mun_\l \ge \eta \frac \eps 2 (\b-\a),$$
so that, if we choose $\eta =\frac \eps 4 (\b-\a)$, we get that
$$ \int f d\nu - \int f d\mun_\l \ge \left[ \frac \eps 4 (\b-\a)\right]^2.$$
And we conclude that, if we take $\d < \left[\eps 4 (\b-\a)\right]^2$, the set $\{\l:d(\mun_\l,\nu)<\d\}$ is empty, which gives (\ref{horsL}).


On the other side, by lemma \ref{defI}.4, for any $M\in\RR^+$,
$$ \tilde\Pi^N(
B(\nu,\d)) 
\le \hspace{-0.3cm}
\sum_{\lambda :d(\mun_\l,\nu)<\d} \hspace{-0.3cm}
I\left(A_N,\phi_M\left({l\over N}\right)\right)^aI\left(B_N,\phi_M\left({l\over N}\right)\right)^b
\D\left({l\over N}\right)^{a+b}  e^{-N^2\int c(x)d\mun_\l(x)
+N^2 F(\mun_\l)}$$
Observe that with $g$ defined in (\ref{defg}),
since $|\l|=\sum\l_j=\sum l_j-\sum(N-j)=\sum l_j -2^{-1}N(N-1)$,
$$ \D\left({l\over N}\right)^{a+b}  e^{-N^2\int c(x)d\mun_\l(x)}
=e^{-N^2\int_{y'\neq y}  g(y',y)
d\mun_\l(y)d\mun_\l(y')
 -N \int c(x)d\mun_\l(x) },$$
we obtain
\begin{eqnarray}
\tilde\Pi^N(
B(\nu,\d))
&\le &\hspace{-0.6cm}
\sum_{\lambda :d(\mun_\l,\nu)<\d} \hspace{-0.3cm}
I\left(A_N,\phi_M\left({l\over N}\right)\right)^aI\left(B_N,\phi_M\left({l\over N}\right)\right)^b\nonumber\\
&&\qquad\qquad \ts \,\, e^{-N^2 \int_{y'\neq y}  g(y',y)\wedge M
d\mun_\l(y)d\mun_\l(y') +N^2F(\mun_\l)
-N \int c(x)d\mun_\l(x)}
\nonumber\\
&\le& e^{N M}
\sum_{\lambda :d(\mun_\l,\nu)<\d} \hspace{-0.3cm}
I\left(A_N,\phi_M\left({l\over N}\right)\right)^aI\left(B_N,\phi_M\left({l\over N}\right)\right)^b
 \nonumber \\
&&\qquad\qquad \ts \,\,e^{-N^2 \int  g(y',y)\wedge M
d\mun_\l(y)d\mun_\l(y') +N^2F(\mun_\l)
-N \int c(x)d\mun_\l(x)}
\label{cutoff}
\end{eqnarray}
Now, following section \ref{IAgrf},
we know that all the functions appearing above are continuous
for any finite $M$ so that for each such $M$  we find a $\kappa(\d,M)$
going to zero as $\d$ goes to zero
 so that

\begin{eqnarray}
\tilde\Pi^N(
B(\nu,\d)) 
&\le& e^{-N^2(\tilde H^M(\nu)+\kappa(\d,M))} e^{N(M+C)}
\sum_{\lambda :d(\mun_\l,\nu)<\d}  e^{-N\rho \int y
d\mun_\l(y) } \label{ml}
\\
\nonumber
\end{eqnarray}
where we used again (\ref{borg}).
We now show that the  last entropy term will not contribute
in the scale $N^2$. We have indeed,

\begin{lem}
\label{tail}
$$ \frac{1}{N^2} \log \sharp\{ \l / d(\mun_\l, \nu) < \delta \} \rightarrow_{N \rightarrow \infty} 0 ,$$
\end{lem}

By (\ref{ml}), and 
lemma \ref{tail} we conclude that, for all $M\ge 0$,

\begin{eqnarray}
\limsup_{N\ra\infty}
{1\over N^2}\log \tilde\Pi^N(
B(\nu,\d))
&\le& -\tilde H^M(\nu) +\kappa(\d,M).
\nonumber
\end{eqnarray}
Letting $\d$ going to zero and then  $M$ going to infinity
(since we saw in section  \ref{IAgrf}
that $\tilde H^M$ converges towards $\tilde H$)
finishes the proof.

\hfill\xx


We now go back to the 
\bf proof of lemma \ref{tail}:\\
\rm We first show a lower bound for the number of tableaux $\l$ whose
 empirical measure is such that, for a given $\epsilon >0$ and a
 given $\nu \in \Ppos$, $d(\frac 1 N \sum_{j=1}^N 
\delta_{\frac{l_j}{N}}, \nu) < \epsilon$.\\
As this number is an integer, we just need to show 
that this set is non-empty. This is true thanks to two facts :
 first the set $\{ \nu \}$ is tight so that we choose a convex compact 
$K$ such that  $\nu(K) \pgq 1 - \frac{\epsilon}{3}$ and then the
 set $\mathcal{P}(K)$ of all probability measures on $K$ endowed 
with the weak topology is a compact in the locally convex space of
 measures with mass less than $1$, so that the Krein-Milman theorem 
tells us that $\mathcal{P}(K)$ is the closure of the convex envelope 
of its extremal points, which are the Dirac measures. We have
 the approximation announced above : 
for $\epsilon > 0$, there exists an integer $N(\epsilon)$ 
and some real number that we order $a_{1, N(\epsilon)} > 
a_{2, N(\epsilon)}> \ldots $ 
such that $d(\frac 1 N \sum_{j=1}^N \delta_{a_{j,N}}, \nu)
 < \frac{\epsilon}{2}$. Then for each $j$ between $1$ and 
$N$, we choose for $l_j$ the integer for which $\frac{l_j}{N}$ 
is the closest from $a_{j,N}$. This gives us that, for $N$ 
large enough
$$  \sharp \left\{ \l/ d\left( \frac 1 N \sum_{j=1}^N 
\delta_{\frac{l_j}{N}}, \nu \right) < \epsilon \right\} \pgq 1.$$
For the upper bound, we first find a compactly supported
 measure $\nu^{\prime}$ (with support $K = [0, M]$) such 
that $d(\nu, \nu^{\prime}) < \frac{\epsilon}{2} $. This 
gives us that
$$ \{ \l/ d(\mun_\l, \nu) < \epsilon \} \subset 
\left\{ \l/ d(\mun_\l, \nu^{\prime}) < 3\frac{\epsilon}{2} \right\}.$$
Let us consider the function $f_2$ given by
\[
f_2(x)  = \left\{ \begin{array}{ll}
 0, & \mbox{if $x \ppq M$} \\
 x-M, & \mbox{if $M \ppq x \ppq M+2N \epsilon$} \\
 2N \epsilon & \mbox{if $x \ge M+2N \epsilon$.}
\end{array}
\right.
\]
$f_2$ is a bounded Lipschitz function whose Lipschitz norm 
is bounded by 1 and such that $\int f_2 d\nu^{\prime} =0$. But, if there exists an $l_j$ greater or equal $2N^2 \epsilon + NM$ 
then $\frac 1 N \sum_{i=1}^{N} f_2\left( \frac{l_i}{N}\right) 
\pgq 2\epsilon   \pgq 3 \frac{\epsilon}{2}$, so that we have
 the inclusion
$$ \{ \l/ d(\mun_\l, \nu) < \epsilon \} \subset \{ \l/ \forall j, \,
 l_j \ppq 2N^2 \epsilon + NM \} $$
and we get the upper bound as we know that
$$ \sharp \{ \l/ \forall j, \, l_j \ppq 2N^2 \epsilon + 
NM \} \ppq (2N^2 \epsilon + NM )^N.$$
Upper and lower bound together give the result announced 
in lemma \ref{tail}. \\

\hfill\xx

\subsection{$(\tilde\Pi^N)_{N\ge 0}$ satisfies 
a large deviation lower bound}\label{ldlb}
In this part we show that
\begin{lem}
$\tilde\Pi^N$ satisfies a large deviation
lower bound,
i.e for any 
 $\nu\in \Ppos$,
$$\liminf_{\d\ra 0}
\liminf_{N\ra\infty} {1\over N^2}\log \tilde\Pi^N(
B(\nu,\d))
\ge -\tilde H(\nu)
.$$
\end{lem}
\prf \
To prove this lower bound, we follow
\cite{BA-G} and consider
discrete approximations of the probability measures 
$\nu\in \{\tilde H<\infty\}$ as follows.
First note that $\tilde H<\infty$ implies that for any $\a < \b$, $\nu([\a,\b]) \le (\b-\a)$.\\
Recall that we saw
at the end of Lemma \ref{exptightlem}
that
 $\tilde H(\nu)\le M$
 implies that for some universal constant
$C$ and $\rho>0$, 
\begin{equation}\label{pl}
\rho \int xd\nu(x)\le M+C \mbox{ and } \,\,
\Sigma(\nu)>-\infty.\end{equation}
The last  condition in particular implies that
$\nu$ have no atoms.
We now construct the following approximations. \\
If $\nu^L=\phi_L\sharp \nu$, by Chebychev inequality,
$$d(\nu, \nu^L)\le \int_{x>L} d\nu\le 
 \rho^{-1}L^{-1}(M+C),$$
and if $\nu$ is in $\mathcal L$, so is $\nu^L$.\\
We then consider 
\begin{eqnarray*}
a_{N,N} & = 
& \inf \left\{ x \,\,/ \,\, \nu^L([0, x]) 
\pgq \frac{1}{N}\right\} \\
a_{i-1,N} & = & \left\{ \begin{array}{ll}
\inf 
\left\{ x \pgq a_{i,N} \,\, / \,\,
 \nu^L((a_{i,N}, x]) \pgq \frac{1}{N} \right\},& \mbox{if $a_{i,N}<L$} \\
 L+ \frac{1}{N},  & \mbox{otherwise.}
\end{array}
\right.
\end{eqnarray*}
It is easy to check that since $\nu$ has no atoms,
for $N\ge N(\eta)$,
\begin{equation}\label{appp2}
d\left(\nu, {1\over N}\sum_{i=1}^N\delta_{a^{i,N}}\right)<\eta+
 \rho^{-1}L^{-1}(M+C).
\end{equation}
Now, for $N,L$ large enough
so that the right hand sides of (\ref{appp2}) 
is smaller that $2^{-1}\d$,
$$\bigcap_{i=1}^N \left\{\left|{l_i\over N}-
a_{i,N}\right|<{\d\over 2}\right\}\subset
 \left\{d\left(\mun_\l,{1\over N}\sum_{i=1}^N\delta_{a^{i,N}}\right)<{\d\over 2}\right\}
\subset \left\{d\left(\mun_\l,\nu\right)<\d\right\}
$$
Therefore
\begin{eqnarray}
\tilde\Pi^N(B(\nu,\d)
&\ge& \tilde\Pi^N\left(\bigcap_{i=1}^N \left\{\left|{l_i\over N}-
a_{i,N}\right|<{\d\over 2}\right\}\right)\nonumber\\
&\ge&\tilde\Pi^N\left(\bigcap_{i=1}^N \left\{\left|{l_i\over N}-
a_{i,N}\right|<\e\right\}\right)\nonumber\\
\nonumber
\end{eqnarray}
for any $\e\in (0,{\d\over 2}]$.
We now show that for any fixed $L$,
\begin{equation}
\label{bi1}
\liminf_{\e\downarrow 0}
\liminf_{N\ra\infty}{1\over N^2}\log
\tilde\Pi^N\left(\bigcap_{i=1}^N \left\{\left|{l_i\over N}-
a_{i,N}\right|<\e\right\}\right) 
\ge -\tilde H(\nu^L)
.  
\end{equation}
Observe first that
${1\over N}\sum_{i=1}^N\delta_{a^{i,N}}$
is supported in $[-L-1,L+1]$
so that all the spherical integrals
are well defined and
uniformly continuous
by Lemma \ref{defI}.
Therefore, we find a $\kappa(\e)$, going to
zero
with $\e$ such that for $N$ sufficiently large,

\begin{eqnarray}
\tilde\Pi^N\left(\bigcap_{i=1}^N \left\{\left|{l_i\over N}-
a_{i,N}\right|<\e\right\}\right)
&\ge&  e^{N^2(aI(\log \sharp\mu_A, \nu^L)
+bI(\log\sharp\mu_B, \nu^L)+F(\nu)-\kappa(\e))} \nonumber\\
&&\ts
\sum_{\left|{l_i\over N}-
a_{i,N}\right|<\e} \D\left({l\over N}\right)^{a+b}
e^{-N^2\int c(x)d\mun_\l(x)
 } \label{enplus}
\end{eqnarray}
Notice that 
\begin{eqnarray*}
\sum_{\left|{l_i\over N}-
a_{i,N}\right|<\e} \hspace{-0.5cm}\D\left({l\over N}\right)^{a+b} e^{-N^2\int c(x)d\mun_\l(x)
 } 
 & =& \hspace{-0.3cm}\sum_{|{l_i\over N}-
a_{i,N}|<\e} \hspace{-0.5cm}
 e^{ N^2 \left( {a+b\over 2}\int \! \!\! \int_{x\neq y} \log |x-y| d\mun_\l
(x)d\mun_\l(y) - \int c(x)d\mun_\l(x)
    \right)} \\
  & \pgq & e^{ -N
\sum_{j=1}^N \sup_{|x-a_{j,N}|\le \frac
{\d}{2}} c(x) 
+{a+b\over 2} N^2 
\int \! \! \int_{x\neq y} \log |x-y| d\mun_\l(x)d\mun_\l(y) }. 
\end{eqnarray*}
where $\l$ is a Young shape defined by
$l_i := \lfloor Na_{i,N} \rfloor$.\\
Note that such a tableau exists since
 according to the definition of the $a_{i,N}$'s since  we have that 
$$ \frac 1 N \le \nu^L([a_{i+1,N},a_{i,N}]) \le a_{i,N} - a_{i+1,N},$$
so that
$$ N(a_{i,N} - a_{i+1,N}) \ge 1,$$
which insures that $l_i-l_{i+1}\ge 1$
and so $\l_i\ge \l_{i-1}$ for all $i\in\NN$.
Note that $|{l_i\over N}-
a_{i,N}|< \frac 1 N$ is smaller than $\eps$ for $N$ large enough.\\
Furthermore, we also get the estimate
 $$ a_{i+1,N} \le \frac{l_i}{N} \le a_{i,N}.$$
Therefore, for $i$, $j$ such that $i<j-1$, we have the lower bound
$$ \left| \frac{l_i}{N} - \frac{l_j}{N} \right| \ge \left|a_{i,N} - a_{j-1,N}\right|,$$
so that we get
\begin{multline*}
\sum_{|{l_i\over N}-
a_{i,N}|<\e} \hspace{-0.5cm}\D\left({l\over N}\right)^{a+b} e^{-N^2\int c(x)d\mun_\l(x)
}
   \pgq  \exp \left( N^2 \left(
 -\frac 1 N \sum_{j=1}^N (c(a_{j,N})+C(L, \d)) 
\right. \right. \\
 +{a+b\over 2} \left. \left. \frac{1}{N^2} 
\sum_{i+1<j} \log |a_{i,N}-a_{j,N}| +
 \frac{a+b}{4 N^2} \sum_{i=1}^{N-1} 
\log \left| \frac{l_{i+1}}{N}- \frac{l_{i}}{N}\right|\right)\right)
\end{multline*}
where we $C(L,\d)$ 
is going to zero as $\d$ goes to infinity
for any given $L$.
With our choice of the $a_{j,N}$'s, we have that
$$ \lim_{N\ra\infty}
\frac{1}{ N} \sum_{j=1}^N c(a_{j,N}) = \int x d\nu^L(x),$$
and 
\begin{eqnarray}
\label{pol}
\lefteqn{\frac{1}{N^2} \sum_{i<j} \log|a_{i,N}-a_{j+1,N}| + \frac{1}{2N^2} \sum_{i=1}^{N-1} \log|a_{i,N}-a_{i+1,N}|} \nonumber \\
& = & \sum_{1\ppq i \ppq j \ppq N-1} 
\log|a_{i,N}-a_{j+1,N}| \nu^L
\otimes \nu^L(a_{i,N} \ppq x \ppq a_{i+1,N};\, a_{j,N} 
\ppq y \ppq a_{j+1,N}) \nonumber \\
& \pgq & \int_{a_{1,N} \ppq x< y \ppq a_{N,N}} \log
|x-y|d\nu^L(x)d\nu^L(y)
\end{eqnarray}
 Let's turn our attention to the last term : 
 for any choice of the $l_i$'s, as the $l_i$ are distinct integers, the difference of a pair of them is at least $1$, so that
 we have
 $$ \prod_{i=1}^{N-1} \left|  
\frac{l_{i+1}}{N}- 
 \frac{l_i}{N}
 \right|  \pgq \left(\frac{1}{N}\right)^{N-1},$$
 which gives     
 $$ \liminf_{N \rightarrow \infty} 
\frac{1}{N^2} \log  \sum_{i=1}^{N-1}  
\log \left|   \frac{l_{i+1}}{N} - \frac{l_i}{N}
 \right| = 0.$$\\
Putting everything together, we can conclude, 
$$
\liminf_{\e\downarrow 0}
\liminf_{N\ra\infty } {1\over N^2}\log 
\sum_{|{l_i\over N}-
a_{i,N}|<\e} \D\left({l\over N}\right)^{a+b} e^{-N^2\int c(x)d\mun_\l(x)}
\ge -{a+b\over 2}\Si(\nu^L)-\int c(x) d\nu^L(x)
$$
(\ref{pol}) and (\ref{enplus}) prove ({\ref{bi1}).
To finish the proof ,
we take the supremum over
$L$ to obtain the lower bound thanks to Lemma \ref{defI}.2 and monotone convergence theorem.\\

\section{Laplace method for $Z_N(\Phi)(X)$ } \label{lap}

Let $\mu^N_\phi$ be the measure on $\Pa(\RR)$ given,
for any measurable set $X$ of $\Pa(\RR)$, 
by
$$\mu^N_\phi(X)={Z_N(\Phi)(X)\over Z_N(\Phi)}.$$
The goal of this section
is to prove a large deviation theorem for $\mu^N_\phi$.\\
We first need some definitions.

\begin{defi}\label{defrf}
With 
$ \mathcal L$ as defined in Theorem \ref{scc} and $\rho_\Phi$
given by (\ref{borne}),
 we let 
$$
\mathcal G_{\Phi}(\nu)  =  \left\{ \begin{array}{ll}
 -I(\log \sharp \mu_A, \nu) - I(\log \sharp \mu_B, \nu) -\Si(\nu)
 +\rho_\Phi. \int xd\nu(x), & 
\mbox{if  $\nu \in \mathcal L$}, \\
+ \infty &  \mbox{otherwise,} 
\end{array} \right.
$$
and if $\Psi=||\Phi||_\infty^{-1}\Phi$,
$$J_\Phi(\nu,\mu):=
   \left\{ \begin{array}{ll}
-I(\log\Psi \sharp \mu,\nu)
-{1\over 2}S(\Psi\sharp\mu) -\Si(\mu)+{1\over 2}\int x^2 d\mu(x), & 
\mbox{if   $\nu \in \mathcal L$}, \\
+ \infty &  \mbox{otherwise.} 
\end{array} \right.
$$
The rate function 
governing our large deviation principle
is then  given, for $\mu\in\Pa(\RR)$, by

$$I_\Phi(\mu):= 
\inf_{\nu \in \Ppos} 
\left( \Ga_{\Phi}(\nu) + J_{{\Phi}}(\nu,\mu) \right)
-\inf_{\mu' \in \Pa(\RR)}\inf_{\nu' \in \Ppos} 
\left( \Ga_{\Phi}(\nu') + J_{{\Phi}}(\nu',\mu') \right)
.$$
\end{defi}

To prove the large deviation principle, we shall 
make the following additional hypothesis
\begin{hyp}\label{hypothese2}
\mbox{}\\
The cut-off function $\Phi$ is bounded below :
\begin{equation}
\label{phi:eps}
\exists \, \epsilon > 0 \, \, s.t. \,\, \forall \,
 x \in \RR, \, \, \Phi(x) \ge \e.
\end{equation}
The two sequences of matrices 
$(A_N)_{N \in \NN}$ and $(B_N)_{N \in \NN}$ and their spectral measures $\hat\mu_{A_N}$ and $\hat\mu_{B_N}$
are such that \\
$\bullet$ there exists an $\a>0$ so that for all $N$, 
$A_N$ and $B_N$ are bounded below by $\a I$.
Hence, with $\Ka$
the compact set $[\a,1]$,
$\textrm{supp} \,\hat\mu_{A_N} \subset \mathcal K$ and  
$\textrm{supp}\, \hat\mu_{B_N} \subset \mathcal K$.  \\
$\bullet$ $\mu_{A_N}$ and $\mu_{B_N}$ 
converge weakly respectively to $\mu_A$ and $\mu_B$.
\end{hyp}
We shall then prove that

\begin{theo}\label{ldp}
Under Hypotheses \ref{hypothese} and  \ref{hypothese2},
\begin{enumerate}
\item $I_\Phi$ is a good rate function on $\Pa(\RR)$,
i.e. $I_\Phi$  is non-negative and
for any $M\in\RR^+$, $\{\nu\in\Pa(\RR): I_\Phi(\nu)\le M\}$
is compact.

\item $(\mu^N_\Phi)_{N\in\NN}$ satisfies a large
deviation principle in the scale $N^2$ with good rate function
$I_\Phi$, i.e

$\bullet$ For any closed subset $F$ of $\Pa(\RR)$,
$$\limsup_{N\ra\infty}{1\over N^2}\log
\mu^N_\Phi(F)\le -\inf_F I_\Phi,$$
$\bullet$ For any open subset $O$ of $\Pa(\RR)$,
$$\liminf_{N\ra\infty}{1\over N^2}\log
\mu^N_\Phi(O)\ge -\inf_O I_\Phi.$$
\item Under Hypothesis \ref{hypothese2}, $S(\mun_{A_N})$ 
converges towards $S(\mu_A)$
and idem for $B_N$, and 

$$\lim_{N \rightarrow \infty}\frac{1}{N^2} \log {Z_N(\Phi)\over
Z_N(0)}
= - \inf_{\mu \in \PR}
\inf_{\nu \in \Ppos} \left( \Ga_{\Phi}
(\nu) + J_{\Phi}(\nu,\mu)  \right)
+{1\over 2} S(\mu_A)+{1\over 2}S(\mu_B)+\frac 1 2 \rho_\Phi.$$
\end{enumerate}

\end{theo}
The proof of this theorem is deduced from
a large deviation principle
obtained for the 
law of the couple $(\mun_\l,\mun_M)$
given by the Gibbs measure defined, for $X=(X_1,X_2)
\subset \Pa(\RR^+)\ts\Pa(\RR)$, by

\begin{equation}
\label{Gibbs}
 \Pi^N_\Phi(X)= {1\over Z_N(\Phi)}
 \sum_{\lambda :\mun_\l\in X_1 } 
s_\l(A_N)s_{\lambda}(B_N)  Z_N(\Psi, \l)(X_2)  e^{-\rho_\Phi |\l|}
\end{equation}
that we can formulate as follows :
\begin{theo}\label{ldp2}
\mbox{}
\begin{enumerate}
\item For $(\nu,\mu)\in \Pa(\RR^+)\ts\Pa(\RR)$, we
set
$$\Ia_\Phi(\nu,\mu):=\left\lbc
\begin{array}{l}
+\infty\mbox{ if } \nu\not\in\La \mbox{ or }  \int
x^2d\mu(x)=+\infty,\\
 J_{\Phi}(\nu,\mu)+\Ga_{\Phi}(\nu) -\inf_{(\nu',\mu')
\in  \Pa(\RR^+)\ts\Pa(\RR)} \{ 
J_{\Phi}(\nu',\mu')+\Ga_{\Phi}(\nu')\}
\mbox{ otherwise.}\\
\end{array}
\right.
$$
Then $\Ia_\Phi$ is a good rate function.
\item $(\Pi^N_\Phi)_{N\in\NN}$
satisfies a full large deviation principle in the
scale $N^2$ 
with rate function $\Ia_\Phi$.

\end{enumerate}
\end{theo}

Theorem \ref{ldp}.1 and  .2  are
direct consequences of Theorem \ref{ldp2} and the contraction 
principle
since 
the application 
$(\nu,\mu)\in  \Pa(\RR^+)\ts\Pa(\RR)\ra \nu\in \Pa(\RR)$
is clearly continuous.

\nn
{\bf Proof of Theorem \ref{ldp2} :}
This proof follows rather closely that
of Theorem \ref{scc2}. Let us briefly outline it.
\begin{enumerate}
\item To prove that $\Ia_\phi$ is a good
rate function, we proceed exactly as in section \ref{grf} ;
$\Ga_\Phi$ has compact level
sets by direct application
of Theorem \ref{scc2}.1 whereas
for $J_\Phi$ we can proceed similarly 
once we notice that $\mu\ra S(\Psi\sharp\mu)$
is continuous since $\Psi$ is bounded below by
a positive constant  and
$$S(\Psi\sharp\mu)=\int\!\!\int \log \left(
\int_0^1(a\Psi(x)+(1-a)\Psi(y))^{-1}da\right)d\mu(x)d\mu(y)$$
and introducing the function
$$j(x,y)=\log |x-y|^{-1}+{1\over 4} x^2+{1\over 4}y^2,$$
we can treat it as $g$ to show that 
$\displaystyle \mu \mapsto \int\!\!\int j(x,y) d\mu(x) d\mu(y)$
is lower semicontinuous on $\PR$.\\ 
Note that we see that $\Ia_\Phi(\nu,\mu)$
is infinite unless
$$\nu\in\La,\quad \int xd\nu(x)<\infty,\quad\Sigma(\nu)>-\infty,
\quad \int x^2 d\mu(x)<\infty, \quad \Sigma(\mu)>-\infty.$$
\item To prove that $\Pi^N_\Phi$ is exponentially
tight, we consider a compact
$$K_L:=\{\nu\in\Ppos : \int xd\nu(x)\le L\} \ts \{\mu\in\Pa(\R):
\int x^2 d\mu(x)\le L\}.$$
It is not hard to bound below $Z^N_\Phi$ by 
some estimate of order $e^{-N^2 C}$ (for instance
by proving the lower bound estimate as below).
Then, using the fact that $S(\Psi\sharp \mu)$ 
is bounded uniformly as well
as the spherical integrals,
we find a finite constant $C'$ such that
$$\Pi^N_\Phi(K_L^c)\le e^{C'N^2}\left(\tilde\Pi^N(\Ka_L^c)+\int_{\sum x_i^2\ge NL}\D(x)^2 e^{-{N\over 2}\sum_{i=1}^N x_i^2}
\prod_{i=1}^N dx_i\right).$$
Following \cite{BA-G} (or the arguments of section \ref{exptight})
we easily see that for sufficiently large $L$
$$\limsup_{N\ra\infty}{1\over N^2}\log \int_{\sum x_i^2\ge NL}\D(x)^2 e^{-{N\over 2}\sum_{i=1}^N x_i^2}
\prod_{i=1}^N dx_i\le -{1\over 4}L$$
so that we can conclude again by section \ref{exptight}.
\item To prove the weak large deviation upper bound, we proceed
as in section \ref{ldub} by considering the functions $g$ (with
$c(x)=\rho_\Phi x$ and $a=b=1$) and $j$.
We then impose a cutoff on both functions and on
the spherical integrals as in (\ref{cutoff})
to obtain a large deviation upper bound estimate, and
then proceed again by optimizing over the cutoff.
\item For the large deviation lower bound, we 
restrict the sum and the integral
also to configurations contained
in small neighborhoods of well 
chosen values $(a_{i,N})_{1\le i\le N}$
and $(x_{i,N})_{1\le i\le N}$
and show convergence.
This strategy works
as well in the continuous setting
as can be seen in \cite{BA-G}.
\end{enumerate}

\section{Comments on the minimizers of $\Ia_\Phi$} \label{secmin}
In this last section, we wish to
give some weak description 
of the minimizers of $\Ia_\Phi$.
We have not been
able to prove uniqueness of such minimizers.
In \cite{G}, uniqueness of the minimizers
of the rate function was deduced from
convexity
arguments which were
actually lacking for instance for the $q$-Potts 
model. In fact, the spherical integrals
are expressed as the sum of 
a convex complicated function  and the entropies $\Sigma$
which are concave. Hence, if the full
rate function does not contain some
term to kill these $\Sigma$ terms,
the convexity of the
full rate function becomes unclear.
The same phenomenon 
appears here and despite our efforts we could not
overcome this difficulty. 
It is unclear here whether the minimizer
should be unique or not.
We here meet the additional difficulty  that
the formula obtained in \cite{GZ}
for the limit of the spherical integral
concerned the case where both 
probability measures had finite covariance, which is not the case here
(one of the argument has
only a first moment which is finite,
even if the other one is
compactly supported).\\
In this section,
we show that the minimizers
of $\Ia_\Phi$ are compactly supported.
We then characterize 
the minimizers. 

\begin{prop}\label{mini} Assume that $\Sigma(\log \sharp \mu_A)>-\infty,\,\,\Sigma(\log \sharp \mu_B)>-\infty.$ Then
\begin{enumerate}
\item 
There exists a 
real number $M\ge 0$ such that
any minimizer $(\nu,\mu)\in\Ppos\ts\Pa(\RR)$ of $\Ia_\Phi$
satisfies $\mbox{supp}(\nu)\subset [0,M]$.
\item If we additionally
 assume
that there exists $A<B$
in $\RR$ such that for $L$ large enough $\Phi$ satisfies
\begin{equation}\label{cc1}
\max_{|x|\ge L} \Phi(x)\le \inf_{x\in [A,B]}\Phi(x)
\end{equation}
then there exists 
a  real number $M$
such that for any minimizer
 $(\nu,\mu)\in\Ppos\ts\Pa(\RR)$ of $\Ia_\Phi$, $\mu$
satisfies $\mbox{supp}(\mu)\subset [-M,M]$.
\item $\Ia_\Phi$ achieves its minimal value (which
is zero).
Let $(\bar\nu,\bar\mu)$
be a minimizer.
Then 
\begin{itemize}
\item There exists 3 flows $(\rho^i,u^i)_{1\le i\le 3}$
such that

$\bullet$ $\mu^i_t(dx)=\rho^i_t(x)dx $ is a 
probability measure for all $t\in (0,1)$.
$t\in [0,1]\ra\mu^i_t\in\Pa(\RR)$ is continuous.
$$\lim_{t\ra 0}\mu^1_t=\log\sharp\mu_A,\quad \lim_{t\ra 0}\mu^2_t=\log\sharp\mu_B,
\quad \lim_{t\ra 0}\mu^3_t=\log\Psi\sharp \mu,$$
$$\lim_{t\ra 1}\mu^i_t=\nu,\quad 1\le i\le 3.$$

$\bullet$ For $i\in\{1,2,3\}$, $(\rho^i,u^i)$ satisfies the
Euler equation
for isentropic flow described by the
equations, 
for $t\in (0,1)$,  

\begin{eqnarray}
\partial_t \rho^i_t(x)&=&-\partial_x( \rho^i_t(x)u^i_t(x))
\label{eulpr10}\\
\partial_t(\rho^i_t(x)u^i_t(x))&=&-\partial_x\left(\rho^i_t(x)
u^i_t(x)^2-{\pi^2\over 3}\rho^i_t(x)^3\right)\label{eulpr20}\\
\nonumber
\end{eqnarray}
in the sense of distributions that
for all $f\in \Ca_c^{\infty,\infty}({\RR}\ts [0,1])$,
$$\int_0^1\int\partial_t f(t,x)d\mu^i_t(x)
dt +\int_0^1\int\partial_x  f(t,x) u^i_t(x) d\mu^i_t(x)dt=0$$
and, for 
any $f\in\Ca^{\infty,\infty}_c(\Omega)$ with
$\Omega_i:=\{(x,t)\in\RR\ts[0,1] :
\rho^i_t(x)>0\}$,
\begin{equation}\label{distf}\int_0^1\int \left(2u^i_t
(x)\partial_t f(x,t)+\left( u^i_t(x)^2
-\pi^2 \rho^i_t(x)^2\right) \partial_x f(x,t)\right)\rho^i_t(x)
dx dt=0,
\end{equation}
where $\Ca^{\infty,\infty}_c(\Aa)$ is the space of functions which are infinitely differentiable on both variables on the open set $\Aa$ and compactly supported.\\
$(\rho^i,u^i)$ are smooth in the interior of $\Omega_i$,
which guarantees that (\ref{eulpr10}) and (\ref{eulpr20})
hold everywhere in the interior of $\Omega_i$. Moreover,
$\Omega_i$ is bounded in $\RR\ts [0,1]$.

$\bullet$
Let $\bar\rho$ be the density of $\bar\nu$
and $\bar\Omega=\{ x: \bar\rho(x)>0\}$
Then, for any continuously
differentiable test function $\phi$
which is supported in the interior 	
of $\bar\Omega$,
$$\int \left(\rho_\Phi x-{1\over 2} x^2 +\int\log|x-y|
d\bar\nu(y)\right)
\partial_x \phi(x) dx =\sum_{i=1}^3\int \phi(x) u^i_1(x)dx.$$
$\bullet$ For any $\phi\in\Ca^1(\mbox{Im}(\log\Psi)^c\cap
\mbox{supp}(\bar\mu))$,
$$\int \partial_x\phi(x)\left( {1\over 2} x^2 -2\int\log|x-y|d\bar\mu(y)
\right) dx=0$$
 To simplify, we shall assume that
$\log\Psi$ is one to
one from $\RR$ into its image $\mbox{Im}(\log\Psi).$ 
Then, in a very weak sense of distribution, for any
$\phi\in\Ca^1(\mbox{Im}(\log\Psi)\cap\mbox{supp}(\bar\mu))$
$$\hspace{-1.5cm}
\int \partial_x \phi\big(
-{1\over 2}x^2 +{1\over 2}(\log\Psi)^{-1}(x)^2
-2\int \log|(\log\Psi)^{-1}(x)-y|d\bar\mu(y)
$$
$$
+\int \log\left|e^x-\Psi(y)\right|d\bar\mu(y)
\big)dx=-\int \phi(x) u^3_0(x) dx.$$
If $\bar\mu$ has a density with respect to
Lebesgue measure, we obtain the usual
sense of distribution in the interior of
$\mbox{Im}(\log\Psi)\cap\mbox{supp}(\bar\mu)$.

\end{itemize}
\end{enumerate}
\end{prop}
The additional assumption is needed to
be able to use \cite{G} results which required it.

\nn
\prf \ \ 
$\bullet$  We first prove the first point, that is for any minimizer $(\nu,\mu)\in\Ppos\ts\Pa(\RR)$ of $\Ia_\Phi$, $\nu$ is compactly supported.
 In \cite{G}, such a result was obtained by
going back to the matrix
model. We shall
here provide a new proof 
based on the study 
of $\Ia_\Phi$. The only property
of the spherical integral we shall
use is the following :
Let $\nu$ and $\nu^*$ in $\Ppos$ be
such that there exists 
a coupling $\pi\in \Pa(\RR^+\ts\RR^+)$ of $(\nu,\nu^*)$ such that
$\pi(x\in .)=\nu(x\in .), \,\, \pi(y\in.) =\nu^*(y\in .)$,
and 
\begin{equation}\label{ploup}
\pi(x\le y)=1.\end{equation}
Then,
for any $\mu\in \Pa(\RR^-)$ which
is compactly supported,
\begin{equation}\label{ineqmin}
I(\nu^*,\mu)\le 
I(\nu,\mu).
\end{equation}
This is a direct consequence of the
definition of the spherical integral ;
indeed, by the above, we can construct 
discrete approximations $(l_i,1\le i\le N)$
and
$(l_i^*,1\le i\le N)$ such that $N^{-1}\sum_{i=1}^N
\d_{l_i\over N}$ (resp. $N^{-1}\sum_{i=1}^N
\d_{l_i^*\over N}$) converges towards $\nu$ (resp. $\nu^*$)
and $l_i\le l_i^*$.
Therefore, if $N^{-1}\sum_{i=1}^N
\d_{\l_i}$ approximates $\mu$ with $\l_i\le 0$,
it is clear that
$$I_N\left({l_i\over N},{\l_i}\right)\ge I_N\left({l_i^*\over N},{\l_i}\right)$$
yielding (\ref{ineqmin}) at the limit $N\ra\infty$.\\
\\
Let now $(\nu^*,\mu^*)$ be
a minimizer and $\nu$  satisfying
(\ref{ploup}) belonging to $\La$.
By definition,
$$\Ia_\Phi(\nu,\mu^*)\ge \Ia_\Phi(\nu^*,\mu^*),$$
and therefore by (\ref{ineqmin}),
 since $\log\sharp\mu_A$,
$\log\sharp\mu_B$ and $\log\Psi\sharp\mu$ 
are supported in $\RR^-$,
\begin{equation}\label{eqmin}
-\Si(\nu) +\rho_\Phi\int x d\nu(x)
\ge -\Si(\nu^*)+\rho_\Phi \int xd\nu^*(x).
\end{equation}


We shall use 
this inequality for a well chosen $\nu$ which is a modification of $\nu^*$. We construct it as follows :
recall that $\nu^*\in\La$
implies that $\nu^*(dx)=\rho^*(x) dx$ with $\rho^*\le 1$.
We assume that $\nu^*([0,M])<1$ 
and are going to show a contradiction 
for $M$ large enough. Observe that 
$A:=\int_0^3  1_{\{x:\rho^*(x)\le {1\over 2}\}} dx \ge 1$
since $\int_0^\infty \rho^*(x) dx =1$.
Set for $M\ge 3$,

$$\nu= \nu_M=1_{[0,M]} \nu^* +{\a_M\over A}1_{\{\rho^*\le {1\over 2}, x\in
[0,3]\}} dx,$$
with $\a_M=\nu^*([M,\infty[)$.\\

We have on one side that 

\begin{eqnarray*}
-\Si(\nu^*)&=& -\Si(1_{[0,M]}\nu^*)
+2\int_{x<M\atop y>M} \log |x-y|^{-1} d\nu^*(x)  d\nu^*(y) 
+\int_{x>M\atop y>M}\log |x-y|^{-1} d\nu^*(x)  d\nu^*(y)\\
&\ge& -\Si(1_{[0,M]}\nu^*)
+2\int_{x<M, y>M\atop |x-y|>1 } \log |x-y|^{-1} d\nu^*(x)  d\nu^*(y) 
\\
&&+\int_{x>M, y>M\atop |x-y|>1}\log |x-y|^{-1} d\nu^*(x)  d\nu^*(y)\\
\end{eqnarray*}
Using that for all $a\in (0,1]$ there exists a
finite constant such that for all $x\ge 0$,
$$\log (1+x)\le C_a x^a$$
we deduce
\begin{eqnarray}
-\Si(\nu^*)&\ge& -\Si(1_{[0,M]}\nu^*)
-2C_a\int_{x<M, y>M\atop |x-y|>1 } (|x-y|-1)^ad\nu^*(x)  d\nu^*(y) 
\nonumber\\
&&-C_a\int_{x>M, y>M\atop |x-y|>1 } (|x-y|-1)^ad\nu^*(x)  d\nu^*(y)
\nonumber\\
&\ge& -\Si(1_{[0,M]}\nu^*) -(2+\a_M)C_a \int_{y>M}
y^a d\nu^*(y) \nonumber\\
&\ge& -\Si(1_{[0,M]}\nu^*) -(2+\a_M)C_a M^{a-1}\int_{y>M} y d\nu^*(y) 
\label{k1}
\end{eqnarray}
\medskip
where we used in the last line Chebyshev 
inequality. \\
On the other side,
\begin{eqnarray}
-\Si(\nu_M)&=& -\Si(1_{[0,M]}\nu^*)
+2{\a_M\over A}
\int_{x<M}\int_0^{3}1_{\rho^*(y)\le {1\over 2}} \log |x-y|^{-1} dy d\nu^*(x)\nonumber\\
&&
+\left({\a_M\over A}\right)^2\int_0^{3}1_{\rho^*(x)\le {1\over 2}}\int_0^{3}
1_{\rho^*(y)\le {1\over 2}}\log |x-y|^{-1} dydx
\nonumber\\
&\le& -\Si(1_{[0,M]}\nu^*)
+2{\a_M\over A}
\int_{x<M}\int_0^{3}1_{\rho^*(y)\le {1\over 2}} 1_{|x-y|\le 1}
\log |x-y|^{-1}dy \rho^*(x) dx\nonumber\\
&&
+\left({\a_M\over A}\right)^2\int_0^{3}1_{\rho^*(x)\le {1\over 2}}\int_0^{3}
1_{\rho^*(y)\le {1\over 2}}1_{|x-y|\le 1} \log |x-y|^{-1} dydx
\nonumber\\
&\le& -\Si(1_{[0,M]}\nu^*)
+\left(2{\a_M\over A}+\left({\a_M\over A}\right)^2\right)
\int_{x<4 } \int_0^3 1_{|x-y|\le 1}\log |x-y|^{-1}  dydx \nonumber\\
&\le& -\Si(1_{[0,M]}\nu^*)
+4\left(2{\a_M\over A}+\left({\a_M\over A}\right)^2\right)\label{b11}\\
\nonumber
\end{eqnarray}

Observe now that $\nu_M$ in $\La$ for $M$ large enough so
that $A^{-1}\a_M \le 2^{-1}.$
Furthermore, $\nu_M$ satisfies (\ref{ploup})
since we have been transporting 
large values of the $l_i$'s
to smaller one.
Hence, we can apply (\ref{eqmin}) and together with (\ref{k1}), (\ref{b11}), it gives that

$$
\rho_\Phi\left(\int_{x>M} xd\nu^*(x)-{\a_M\over A}
\int_0^3 x 1_{\rho^*<{1\over 2}}
dx\right)\le 
(2+\a_M)C_a M^{a-1}\int_{y>M} y d\nu^*(y) +
4\left(2{\a_M\over A}+\left({\a_M\over A}\right)^2\right),
$$
showing that for any $a\in (0,1)$, for $M$ large
enough,
\begin{equation}
(\rho_\Phi -(2+\a_M)C_a M^{a-1})
\int_{x>M} xd\nu^*(x)\le  {15\over A} \a_M\le
{15\over A M} \int_{x>M} xd\nu^*(x)
\end{equation}
which shows that $\int_{x>M} xd\nu^*(x)$
has to be null when $\rho_\Phi -(2+\a_M)C_a M^{a-1}-{15\over A M}>0$
that is for $M$ large enough.\\

$\bullet$ We now pass to the proof of the second point
of the proposition.
Let, with $\b_M=\mu^*([-M,M]^c)$, for $B>A$,
$$\mu_M(dx)=1_{[-M,M]} \mu^*(dx) +{\b_M\over B-A} 1_{[A,B]} dx$$
Because of our assumption, we see
that if $M$ is large enough and $[A,B]$ chosen
so that
$$\inf_{[A,B]}\Phi\ge \sup_{[-M,M]^c}\Phi$$
for any $\nu\in\Ppos$,
$$I(\log\Psi\sharp\mu_M,\nu)\ge I(\log\Psi\sharp\mu^*,\nu).$$
Hence, 
when $(\mu^*,\nu^*)$
minimize $\Ia_\Phi$,
we obtain
\begin{equation}\label{eqmin2}
-\Si(\mu^*)+{1\over 2}\int x^2d\mu^*(x)
-{1\over 2}S(\Psi\sharp\mu^*)
\le -\Si(\mu_M)+{1\over 2}\int x^2d\mu_M(x)
-{1\over 2}S(\Psi\sharp\mu_M)
\end{equation}
Arguing as above,
we find
that, for any $a\in (0,2)$, there
exists a finite constant $C_a$ such that
\begin{eqnarray}\Si(\mu^*)-\Si(\mu_M)&\le &C_a M^{a-2}
\int x^2 d\mu^*(x)\\
-S(\Psi\sharp\mu_M)+S(\Psi\sharp\mu^*)&\le &C\b_M
\label{c2}\\
\nonumber
\end{eqnarray}
where we observed in the last line that 
$\Psi$ was bounded uniformly above and below.
Hence,
we arrive at
\begin{equation}\label{eqmin3}
\left({1\over 2}-C_a M^{a-2}\right)\int_{x\ge M} x^2 d\mu^*(x)
\le C'\b_M\le C'M^{-2}\int_{x\ge M} x^2 d\mu^*(x)
\end{equation}
where $C'=C+B^2$. This is again a contradiction
for sufficiently
large $M$.\\

$\bullet$ We finally study the characterization
of the minimizers. In \cite{G}, the characterization
was done by going back to
the matrix model description.
We shall here  tackle this
problem by a direct study
of the rate function.
Note that by point 1., any minimizers $(\bar\nu,\bar\mu)$
is such that  $\bar\nu$ is compactly supported.
Moreover $\log\Psi\sharp \bar\mu, \log\sharp\mu_A$ and
$\log\sharp\mu_B$
are also compactly supported by our hypotheses  
so that we can apply Property 2.2 in
\cite{G} which says that
 if $\mu,\nu$ are two
probability measures with finite covariance and such that
$\Si(\mu)>-\infty$, $\Si(\nu)>-\infty$, 
\begin{equation}\label{G1}
I(\mu,\nu)=-{1\over 2}\inf_{(\rho,u)\in C(\mu,\nu)}
\{ S(\rho,u)\}
-{1\over 2}\left(\Si(\mu)+\Si(\nu)-\mu(x^2)-\nu(x^2)\right)+c
\end{equation}
where 
$$S(\rho,u):= \int_0^1 \int u_t(x)^2 \rho_t(x) dx dt 
+{\pi^2\over 3} \int_0^1 \int \rho_t(x)^3 dx dt,$$

\begin{multline*}
C(\mu,\nu)=\Bigg\{ \rho_.\in L^1(dxdt),
\int \rho_t(x) dx=1 \,\, \forall t\in [0,1],
\lim_{t\ra 0} \rho_t(x) dx=\mu,
\lim_{t\ra 1} \rho_t(x) dx=\nu, \\
  \partial_t\rho_t(x)+\partial_x(\rho_t(x)
u_t(x)) =0\Bigg\},
\end{multline*}
where the last equality is to be understood
in the sense of distributions. It was shown in 
\cite{G} that the infimum defining
$I$ is achieved at a unique $(u^*,\rho^*)\in C(\mu,\nu)$
which is described by an isentropic
Euler equation with negative pressure $p(\rho)=-{\pi^2\over 3}\rho^3$.
$c$ is a universal constant.
As a consequence of this
formula, since $\Ia_\Phi(\mu,\nu)<\infty$
implies that $\Si(\mu)>-\infty$, $\Si(\nu)>-\infty$
and $\mu(x^2)<\infty$, for any $\nu\in\Ppos$
such that $\nu(x^2)<\infty$,
we find
that

\begin{eqnarray}
\Ia_\Phi(\mu,\nu)&=&
\inf_{((\rho^i,u^i)\in C(\mu^i,\nu))_{ 1\le i\le 3}}
\left\{{1\over 2} \sum_{i=1}^3 S(\rho^i,u^i)
+{1\over 2}\Si(\nu)-\Si(\mu)+{1\over 2}\Si(\Psi\sharp\mu) \right.\nonumber\\
&&
\left.+{1\over 2}\mu(-\log\Psi(x)^2+x^2)-{3\over
2}\nu(x^2)+\rho_\Phi\nu(x)+K(\mu_A,\mu_B)\right\}\nonumber\\
&:=& \inf_{((\rho^i,u^i)\in C(\mu^i,\nu))_{ 1\le i\le 3}}
\Xi\left( ( \rho^i,u^i)_{1\le i\le 3},\nu,\mu \right),
\label{G2}\\
\nonumber
\end{eqnarray}
where $\mu^1=\log \sharp\mu_A$, $\mu^2=\log \sharp\mu_B$,
$\mu^3=\log\Psi\sharp\bar\mu$ and $K(\mu_A,\mu_B)$ is a constant depending only on $\mu_A$ and $\mu_B$ .\\

We now  consider a minimizer $((\bar \rho^i,\bar u^i)_{1\le i\le 3},
\bar\mu,\bar\nu)$ of $\Xi$
in $\Omega:=\{ \nu\in\La, \mu\in\Pa(\RR),(\bar \rho^i,\bar u^i)_{1\le i\le 3}
\in 
C(\log \sharp\mu_A,\nu)\ts C(\log \sharp\mu_B,\nu)\ts C(\log\Psi\sharp\mu,\nu)
\}$. 
To characterize this minimizer,
we perform a small perturbation.
Let $((\rho^i_\e,u^i_\e)_{1\le i\le 3},
\mu_\e,\nu_\e)\in\Omega$ be given, for  compactly supported 
 functions
$(\phi^i)_{1\le i\le 3}$ in $\Ca^{1,1}(\RR\ts [0,1])$
by
$$\rho^i_\e(t,x)=\bar\rho^i(t,x)+
\e\partial_x\phi^i(t,x) \textrm{\ and \ } u^i_\e(t,x)\rho^i_\e(t,x)
=\bar u^i(t,x)\bar\rho^i(t,x)-\e\partial_t\phi^i(t,x),$$
with $\partial_x\phi^i(1,x)=\partial_x\phi(1,x)$
independent of $i$, $\partial_x\phi^i(0,x)=0$
for $i=1,2$.\\
Note that, once we chose the perturbation for $\rho^i$, the form of the perturbation for $u^i\rho^i$
taken above ensures that the first equation $\partial_t\rho^i(t,x) = - \partial_x(u^i(t,x)\rho^i(t,x))$
is automatically satisfied.\\
This implies also
$$\nu_\e=\bar\nu+ \e \partial_x \phi^i(1,x) dx$$
and $$\log\Psi\sharp \mu_\e(dx)=\log\Psi\sharp \bar
\mu(dx)+ \e \partial_x \phi^3(0,x) dx.$$
We perturb
more generally $\mu$
by setting
$$\mu_\e(dx)=\bar\mu(dx)+\e\partial_x\psi(x)dx$$
with the condition
$$\int f(\log\Psi(x)) \partial_x\psi(x)dx
=\int f(x) \partial_x \phi^3(0,x) dx$$
for all bounded continuous
functions $f$.\\
 
We shall assume that 
$$L(\phi) =\sum_{i=1}^3 \int_0^1\int  \left({|\partial_t \phi^i(t,x)|^2
\over \bar\rho^i(t,x)}\right)dxdt
+\sum_{i=1}^3\sup_{t\in (0,1)}  
\left\|{\partial_x \phi^i(t,x)\over\bar\rho^i(t,x)} \right\|_\infty <\infty.$$
It is not hard to
see that under such conditions,
$\Xi((\rho^i_\e,u^i_\e)_{1\le i\le 3},
\mu_\e,\nu_\e)$ is finite.\\
By the condition
$$\Xi((\rho^i_\e,u^i_\e)_{1\le i\le 3},
\mu_\e,\nu_\e)\ge \Xi((\bar \rho^i,\bar u^i)_{1\le i\le 3},
\bar\mu,\bar\nu)$$
we obtain, taking the limit $\e\ra 0$, that

\begin{eqnarray}
&&\int \left(\rho_\Phi x -{3\over 2}x^2\right)  \partial_x \phi(1,x) dx
-{1\over 2}\int x^2 \partial_x \phi^3(0,x) dx
+{1\over 2}\int x^2 \partial_x\psi(x)dx
\nonumber\\
&&+ \int\!\!\int\log|x-y|d\bar\nu(y)\partial_x\phi(1,x)
dx- 2\int\!\!\int\log|x-y|d\bar\mu(y) 
\partial_x\psi(x)dx\nonumber\\
&&+\int\!\!\int\log\left|e^x-e^y\right|d\log \Psi \sharp\bar\mu(y) 
\partial_x\phi^3(0,x)dx\nonumber\\
&&+{1\over 2}\sum_{i=1}^3 
\int \int_0^1 [-2\partial_t \phi^i(t,x) \bar u^{i}(t,x) 
 -(\bar u^i(t,x) )^2\partial_x \phi^i(t,x) 
+\pi^2(\bar\rho^i(t,x) )^2\partial_x\phi^i(t,x) ]
dxdt\ge 0 \nonumber \\
&&
\label{eqimin}\\
\nonumber
\end{eqnarray}
Changing $\phi^i$ (for $1\le i \le 3$) and $\psi$ respectively into  $-\phi^i$ and $-\psi$, we get that the inequality in (\ref{eqimin}) is in fact an equality.\\
Applying this result with $\phi^i(0,x)=\phi^i(1,x)=0$
shows that $(\bar u^i,\bar \rho^i)_{1\le i\le 3}$
satisfies the Euler equation
for isentropic flow
described in the proposition.\\
\\

 We now turn to the boundary conditions expressed in the last two points of Proposition \ref{mini}.
 To
characterize them, we will try to regularize the densities $\rho^i_\e(t,.)$.
We remark that 
by Property 2.8 in \cite{G},
since $\bar\nu$ and $\bar\mu$ are compactly supported
under our hypothesis, we can find sequences of potentials $
(h^{\e,i}, \e>0, 1\le i\le 3)$
in $\Ca^{1,1}_b(\RR\ts [0,1])$
such that
if we set 
$$\rho^i_\e(t,x):=\pi^{-1} 
( \max\{
\partial_th^{\e,i}(t,x)+4^{-1}(\partial_xh^{\e,i}(t,x))^2,0\} )^{1\over 2}$$
then for any $\e>0$,

$$
 \int \left(\bar u^i(t,x)-\partial_x{h^{\e,i}(t,x)\over 2}
\right)^2 \bar \rho^i(t,x) dxdt+ {\pi^2\over 3}\int_0^1\int \left( \bar\rho^i(t,x)-
\rho_\e^i(t,x)\right)^2 
\left(\bar\rho^i(t,x)+\rho_\e^i(t,x) \right)dxdt$$
$$
+ \pi^2\int_0^1\int |\partial_th^{\e,i}(t,x)
+4^{-1}(\partial_xh^{\e,i}(t,x))^2 -\pi^2
\rho_\e^i(t,x)^2|  \bar \rho^i(t,x) dxdt
 \le\e.$$
From
this result, we deduce
that
$$\sup_{1\le i\le 3}
\left|\int_0^1\int [-2\partial_t \phi^i(t,x) \bar u^{i}(t,x) 
 -(\bar u^i(t,x) )^2\partial_x \phi^i(t,x) 
+\pi^2(\bar\rho^i(t,x) )^2\partial_x\phi^i(t,x) ]
dxdt \right.$$
$$- \left.\int_0^1\!\!\int [-\partial_t \phi^i(t,x) \partial_x h^{\e,i}(t,x) 
 -\frac 1 4 (\partial_x h^{\e,i}(t,x))^2\partial_x \phi^i(t,x) 
+\pi^2(\rho_\e^i(t,x) )^2\partial_x\phi^i(t,x) ]
dxdt\right|\le C_{L(\phi)}\sqrt{\e}$$
with $C(L(\phi))<\infty$ when $L(\phi)<\infty$.
Moreover, since $h^{i,\e}\in\Ca^{1,1}(\RR\ts[0,1])$, we can integrate by part so that
$$\left| \int_0^1\int [-\partial_t \phi^i(t,x) \partial_x h^{\e,i}(t,x) 
 -4^{-1}(\partial_x h^{\e,i}(t,x))^2\partial_x \phi^i(t,x) 
+\pi^2(\rho_\e^i(t,x) )^2\partial_x\phi^i(t,x) ]
dxdt \right.$$
$$-2\left.\left[\int h^{\e,i} \partial_x \phi^i dx\right]_0^1\right|\le
C'(L(\phi))\sqrt{\e}$$
We now can define
in the sense of distribution
$$\int \Pi_t^i \partial_x \phi^i dx
=-\int u_t^i  \phi^i dx$$
and by letting
$\e$ going to zero we get
that
$$\int [-2\partial_t \phi^i(t,x) \bar u^{i}(t,x) 
 -(\bar u^i(t,x) )^2\partial_x \phi^i(t,x) 
+\pi^2(\bar\rho^i(t,x) )^2\partial_x\phi^i(t,x) ]
dxdt=2\left[\int \Pi_t^i \partial_x \phi^i dx\right]_0^1.$$
Thus,
we have proved that we can rewrite (\ref{eqimin}) (which we showed to be an equality) under the form
\begin{eqnarray}
&&\int \left(\rho_\Phi x -{3\over 2}x^2\right)  \partial_x \phi(1,x) dx
-{1\over 2}\int x^2 \partial_x \phi^3(0,x) dx
+{1\over 2}\int x^2 \partial_x\psi(x)dx
\nonumber\\
&&+ \int\!\!\int\log|x-y|d\bar\nu(y)\partial_x\phi(1,x)
dx- 2\int\!\!\int\log|x-y|d\bar\mu(y) 
\partial_x\psi(x)dx\nonumber\\
&&+\int\!\!\int\log|\Psi(x)-\Psi(y)|d\bar\mu(y) 
\partial_x\psi(x)dx\nonumber\\
&&+\sum_{i=1}^3\left(\int \Pi_1^i \partial_x \phi (1,x)dx
-\int \Pi_0^i \partial_x \phi^i (0,x)dx\right)
= 0
\label{eqimin2}\\
\nonumber
\end{eqnarray}
From that we can deduce the
boundary conditions we are seeking for.\\
As the equality (\ref{eqimin2}) holds for any function $\partial_x \phi(1,x)$ such that
$L(\phi)$ is finite,
we find
that 

\begin{equation}\label{dist1}
A(x,\bar\nu)= \rho_\Phi x-{3\over 2} x^2 +\int\log|x-y| d\bar\nu(y)
+\sum_{i=1}^3\Pi^i_1(x)
\end{equation}
is constant in the sense of distribution.\\
Furthermore, it is
not hard to deduce from the representation of $\rho_t^i$
as a free Brownian motion 
given in \cite{G} that for $t$ close enough to
one $\{x: \rho_t^i(x)\ge \e\}\subset \{x: \bar \rho(x)\ge 2\e\}$
with $\bar\rho$
the density of $\bar\nu$
with respect to
Lebesgue
measure. Therefore, 
for any $\Ca^{1}_b$ function
$\phi$ with compact support in
the interior
of $\{x: \bar \rho(x)>0\}$,
$$\int \partial_x \phi(x) A(x,\bar\nu) dx=0.$$

Now only the last point of our proposition is left to establish.\\
The statement of the result 
is more obscur when dealing 
with $\bar\mu$ since we do not a priori
know
if $\bar\mu$ has a density with respect to
Lebesgue
measure. What we get from (\ref{eqimin2})
is that :

 For any $\psi\in
\Ca^1_b(\mbox{Im}(\log\Psi)^c\cap\mbox{supp}(\bar\mu))$
$$\int \partial_x\psi(x)\left({1\over 2}x^2-2\int\log|x-y|d\bar\mu(y)\right) dx=0$$
i.e ${1\over 2}x^2-2\int\log|x-y|d\bar\mu(y)$
is constant outside of the image $\mbox{Im}(\log\Psi)$
 of $\log\Psi$.\\
Inside $\mbox{Im}(\log\Psi)$, if we assume that $\log\Psi$ is one to one from $\RR$ onto its image, we have that
$$B(x,\bar\mu)=-{1\over 2}x^2 +{1\over 2}(\log\Psi)^{-1}(x)^2
-2\int \log|(\log\Psi)^{-1}(x)-y|d\bar\mu(y)
+\int \log\left|e^x-\Psi(y)\right|d\bar\mu(y)
-\Pi^3_0(x)$$
is constant in the weak sense 
of distribution that is its integral with respect
to $\partial_x\phi^3(x,0)$ vanishes. 
If $\bar\mu$ has a
density with respect to
Lebesgue
measure, we find that $B(x,\bar\mu)$
is constant in the sense of
distribution
inside $\{x:{d\bar\mu\over dx}\neq 0\}$
as above, but it is not clear that
a $\phi^3\neq 0$ indeed exists 
in general !

\hfill\xx

\bigskip

\bigskip
\section{Conclusion
and remarks }\label{KSWsec}

In this paper, we studied 
the asymptotics of the 
 model given 
by the partition function (\ref{fnpart1}).
In the course of doing so, 
we adapted the techniques 
of \cite{BA-G}
to study large deviations 
of the profiles of Young tableaux
with a density 
given by a Vandermonde
determinant and Schur  polynomial
functions (see Theorem \ref{scc}). We believe
that these techniques might
be useful to study
other problems
since these kind of distributions
appear in different contexts
due to their  combinatorial
nature. For instance,  following  Migdal-Witten 
formula \cite{Wi,Wi2},
the partition function 
of two-dimensional Yang Mills theory
on a cylinder
with gauge group $U(N)$
is given by 
 the central heat kernel defined, at  time $t=TN^{-1}$,
by
$$\Za_N\left(U_1,U_2; {T\over N}\right)
=\sum_{\l}
s_\l(U_1)s_\l(U_2) e^{-{T\over 2N} C_2(\l)}
$$
where $U_1,U_2\in U(N)$, the sum runs over Young tableaux $\l$ and
$$ C_2(\l)= \sum_{i=1}^N 
\l_i(\l_i+1-2i+N)
=\sum_{i=1}^N l_i^2-(N-1)\sum_{i=1}^N l_i+\sum_{i=1}^N (N-i)(i-1)$$
with $l_i=\l_i+N-i$ (see for example \cite{GM}).\\
S. Zelditch  \cite{zeld} asked us
if we could study the asymptotics of $\Za_N(U_1,U_2; TN^{-1})$
when $U_1,U_2$ are not unitary but
real diagonal
matrices with converging 
spectral distributions.
Our techniques apply readily
to this context and we find 

\begin{theo}
\label{zeld}
Let $A_N,B_N$
be two sequences
of uniformly bounded matrices bounded
below by $\e I$ for some $\e>0$
with spectral measures 
converging towards $\mu_A,\mu_B$.
Then for any time $T>0$
$$\lim_{N\ra\infty}
{1\over N^2}\log \Za_N\left(A_N,B_N; {T\over N}\right)
=Z(\mu_A,\mu_B,T)$$
with 
\begin{multline*}
Z(\mu_A,\mu_B;T)
=\sup_{\nu\in\La}\left\{ I(\log\sharp\mu_A,\nu)
+I(\log\sharp\mu_B, \nu)+\Sigma(\nu)-{T
\over 2}\int x^2 d\nu(x)+\frac T 2 \int xd\nu(x)\right\} \\
+{1\over 2}S(\mu_A)+{1\over 2}S(\mu_B)-{T\over 12}
\end{multline*}
\end{theo}
This theorem is a direct consequence 
of Theorem \ref{scc} with $a=b=1$
and $c(x)=x^2-x$.

\medskip

In addition to giving a rigorous
 basis to the study of 
such natural asymptotics,
we gave a firm ground 
to begin the study
of other matrix models
where other problems
due for instance to signed series
might appear. This step
seems necessary since the proofs are already rather
involved. Furthermore,
 we developed
new arguments to study the saddle points
of our model based on transport of mass.

\medskip

One of the weakness of 
our result is apparently 
the
 cut-off function $\Phi$,
since the matrix integral (\ref{fnpart1})
is then hard to relate with the enumeration of maps
as in \cite{KSW}. Let us comment heuristically
this point. Observe first  that the matrix integral (\ref{fnpart1})
with $\Phi(x)=x$ considered in \cite{KSW}
is always infinite. Indeed, 
for instance in the case $A=1$,
we are integrating 
$$Z_N(Id)=\int_{x_i\in\R} \D(x
)^2\prod_{i,j=1}^N{1\over 1-b_i x_j} e^{-N\sum x_j^2}\prod dx_j$$
which is clearly infinite for all $N\in\NN^*$. Hence,
everything should be understood
formally. The same problem
a priori also arise when one considers
random triangulations generated 
by the one matrix integrals
$$\tilde Z_N(\l)=\int e^{\l N\tr(M^3)-{N\over 2}\tr(M^2)}
dM$$
which is clearly infinite for $\l\neq 0$.
One way to bypass this problem 
is for instance to consider
$$\tilde Z_N(\l,\eta)=\int e^{-\eta N\tr(M^4)+\l N\tr(M^3)
-{N\over 2}\tr(M^2)} dM$$
which is well defined for $\eta>0$.
Recall that planar maps are enumerated by
$$C(n)=\lim_{N\ra \infty}\partial_\l^n
{1\over N^2}\log \tilde Z_N(\l)|_{\l=0}
=\lim_{N\ra \infty}\partial_\l^n
{1\over N^2}\log \tilde Z_N(\l,\eta)|_{\l=0,\eta=0}.$$
In the physics literature, these quantities are implicitely 
supposed to be given by 
$$\tilde C(n)=\partial_\l^n \lim_{N\ra \infty}
{1\over N^2}\log \tilde Z_N(\l,\eta)|_{\l=0,\eta=0}.$$
This 
seems to be fine in the one matrix
case after the work of N.~Ercolani and K.~McLaughlin \cite{EML}
but this point is open in general.

Similarly, one could try to 
regularize the dually weighted graph 
model by considering \linebreak
$Z_N(\Phi_{\e,R})$
with $$\Phi_{\e,R}(x)={x\over 1+\e x^2}+R$$
with $\e>0$ and $R\ge \sqrt{2\e}^{-1}$.
For $||A||$ and $||B||$ small enough (which we can always assume
since again only derivatives at the origin
should be of interest), we obtain by our result
a limit for $N^{-2}\log Z_N(\Phi_{\e,R})$.
Assuming that the limit 
can be extended analytically
to $R,\e$ small,
we should be able to enumerate,
modulo the above ansatz of interchanging 
derivation and limit, 
the enumeration of
dually weighted graphs.

There is still a long way toward
the rigorous understanding of the
use of matrix integrals for the enumeration
of  maps in physics but
we hope that this paper provides
some useful  steps in this direction.

\bigskip

\noindent
{\bf Acknowledgments :}
A.~Guionnet wishes to
thank A.~Okounkov for  patient and cheerful
discussions around  
character expansions.
We also wish
to thank S. Zelditch for 
explaining us the
problem tackled
in section \ref{KSWsec}. 
 M.~Ma\"{\i}da is very grateful to O.~Zeitouni for encouraging and helping her to read papers of physicists and both authors are very indebted to him for various stimulating discussions and always pertinent remarks.
 
\bibliographystyle{plain}
\bibliography{dwg}

\end{document}